\documentclass[12pt,english]{article}
\usepackage[T1]{fontenc}
\usepackage[latin1]{inputenc}
\usepackage{a4}
\usepackage{amsmath}
\usepackage{graphicx}
\usepackage{setspace}
\usepackage{amssymb}

\makeatletter

\newcommand{\noun}[1]{\textsc{#1}}
\providecommand{\tabularnewline}{\\}

 \usepackage{verbatim}

\newcommand{\trace}{\mathop{\mathrm{Tr}}}

\newcommand{\bic}{\mathop{\mathrm{BIC}}}
\usepackage{color}

\usepackage{psfrag} 
\usepackage{epsfig}
\usepackage{natbib}
\usepackage{fancyhdr}

\date{}
\usepackage{babel}
\makeatother
\begin{document}

\title{High-Dimensional Data Clustering}

\author{C. Bouveyron$^{\textrm{a,b,c}}$, S\noun{.} Girard$^{\textrm{a,b}}$
\noun{}and \noun{}C. \noun{}Schmid$^{\textrm{b}}$\\
\vspace{-1ex}{\small $^{\textrm{a}}$}\emph{\small LMC-IMAG, Université Grenoble
1, BP. 53, 38041 Grenoble Cedex 9, France}\\
\vspace{-1ex}{\small $^{\textrm{b}}$}\emph{\small INRIA Rhône-Alpes, 655 avenue
de l'Europe, 38330 Saint-Ismier Cedex, France}\\
\vspace{-1ex}{\small $^{\textrm{c}}$}\emph{\small Dpt. of Mathematics and Statistics, Acadia University, Wolfville, B4P 2R6, Canada}}

\maketitle
\noindent \thispagestyle{fancy}

\lfoot{\textit{Technical report 1083M, LMC-IMAG}}

\cfoot{}\rfoot{\textit{Revised on \today}}

\hrule\medskip\begin{singlespace}\small{

\noindent \textbf{Abstract\bigskip}

\noindent Clustering in high-dimensional spaces is a difficult problem
which is recurrent in many domains, for example in image analysis.
The difficulty is due to the fact that high-dimensional data usually
live in different low-dimensional subspaces hidden in the original
space. This paper presents a family of Gaussian mixture models designed
for high-dimensional data which combine the ideas of subspace
clustering and parsimonious modeling. These models give rise to a clustering
method based on the Expectation-Maximization algorithm which is called
High-Dimensional Data Clustering (HDDC). In order to correctly fit
the data, HDDC estimates the specific subspace and the intrinsic dimension
of each group. Our experiments on artificial and real datasets show
that HDDC outperforms existing methods for clustering high-dimensional
data.\textbf{\bigskip}


\noindent }\end{singlespace}\bigskip\hrule

\section{Introduction}

\noindent Clustering in high-dimensional spaces is a recurrent problem
in many fields of science, for example in image analysis. Indeed,
the data used in image analysis are often high-dimensional and this
penalizes clustering methods. In this paper, we focus on model based
approaches, see~\cite{Bock96} for a review on this topic.
Popular clustering methods are based
on the Gaussian Mixture Model (GMM)~\cite{McLachlan00a} and show
a disappointing behavior when the size of the dataset is too small
compared to the number of parameters to estimate. This well-known
phenomenon is called \emph{curse of dimensionality} and was introduced
by Bellman~\cite{Bellman57}. We refer to~\cite{Pavlenko03,Pavlenko01}
for a theoretical study of the effect of dimension in the supervised
framework. 

\noindent To avoid over fitting, it is necessary to find a balance
between the number of parameters to estimate and the generality of
the model. We propose a Gaussian mixture model which
takes into account the specific subspace around which each cluster is
located and therefore limits the number of parameters to estimate.
The Expectation-Maximization (EM) algorithm~\cite{Dempster77} is
used for parameter estimation and the intrinsic dimension of each
group is determined automatically thanks to the BIC criterion and 
the scree-test of Cattell.
This allows to derive a robust clustering method in high-dimensional
spaces, called High Dimensional Data Clustering (HDDC). In order to
further limit the number of parameters, it is possible to make additional
assumptions on the model. For example, it can be assumed that classes are
spherical in their subspaces or fix some parameters to be common between
classes. The nature of the proposed parametrization allows HDDC 
to be robust with respect to the ill-conditioning or
the singularity of empirical covariance matrices and to be efficient in terms of computing time. Finally, HDDC is evaluated and compared to standard clustering
methods on artificial and real datasets. 

This paper is organized as follows. Section~\ref{sec:Related-work}
presents the state of the art on clustering of high-dimensional data.
Section~\ref{sec:A-Gaussian-mixture} introduces our parameterization
of the Gaussian mixture model. Section~\ref{sec:Clustering-method}
presents the clustering method HDDC, \emph{i.e.} the estimation of
the parameters of the models and of the hyper-parameters. Experimental
results on simulated and real datasets are
reported in Section~\ref{sec:Experimental-results}.

\section{\label{sec:Related-work}Related work on high-dimensional clustering}

Standard methods to overcome the \emph{curse of dimensionality} consist
in reducing the dimension of the data and/or to use a parsimonious
Gaussian mixture model. More recently, methods which find clusters
in different subspaces have been proposed. In this section, 
a brief survey of these works in clustering of high-dimensional data
is presented.

\subsection{Dimension reduction}

Many methods use global dimension reduction techniques to overcome
problems due to high dimensionality. A widely used solution is to
reduce the dimension of data before using a classical clustering
method. Dimension reduction techniques can be divided into techniques
for \emph{feature extraction} and \emph{feature selection}. Feature
extraction techniques build new variables carrying a large part of
the global information. Among these techniques, the most popular one is
Principal Component Analysis (PCA)~\cite{Jolliffe86} which is often
used in data mining and image analysis. However, PCA is a linear technique,
\emph{i.e.} it only takes into account linear dependences between
variables. Recently, many non-linear techniques have been proposed
such as Kernel PCA~\cite{Scholkopf98}, non-linear PCA~\cite{Girard00,Hastie89}
and neural networks based techniques~\cite{Demartines97,Kohonen95,Roweis00,Tenenbaum00}.
In~\cite{Schott93}, the dimension reduction problem was considered
in the Quadratic Discriminant Analysis framework. In contrast, feature
selection techniques find an appropriate subset of the original variables
to represent the data. A survey on feature selection can be found
in~\cite{Guyon03}. A recent approach~\cite{Raftery05} proposes
to combine global feature selection and model-based clustering. These
global dimension reduction techniques are often advantageous in terms
of performance, but suffer from the drawback of losing information
which could be discriminant. Indeed, the clusters are usually hidden
in different subspaces of the original feature space and a global
approach cannot capture this.

\subsection{Parsimonious models}

Another solution is to use models which require the estimation of
fewer parameters. For example, the eigenvalue decomposition of the
covariance matrices~\cite{Banfield93,Celeux95} allows to re-parameterize
the covariance matrix of the classes in their eigenspaces. By fixing
some parameters to be common between classes, this parameterization
yields parsimonious models which generate clustering methods based
on the EM algorithm. A review on parsimonious models can be found
in~\cite{Fraley02}. These approaches are based on various Gaussian
models from the most complex one (a full covariance matrix for each
group) to the simplest one (a spherical covariance matrix for all
groups) which yields a method similar to the k-means
approach. However, these methods
cannot efficiently solve the problem of the high-dimensionality when clusters
live in low-dimensional subspaces.

\subsection{Subspace clustering}

Subspace clustering methods involve two kinds of approaches.
On the one hand, projection pursuit clustering assumes that 
the class centers are located on a same unknown subspace~\cite{Bock87,Soete94}. 
On the other hand, principal component clustering assumes that
each class is located on a unknown specific subspace, see~\cite{Bock74},
Chapter~17, and~\cite{Bezdek81} for an extension to fuzzy subspaces.
For instance, the {\it Analyse factorielle typologique}~\cite{Diday74}
is based on an iterative algorithm similar to the k-means approach.
Some subspace clustering methods use heuristic search techniques
to find the subspaces, see for instance~\cite{Agrawal98}. 
A review on this type of methods can be found in~\cite{Parsons04}.
Most of them rely on geometric considerations and are not model-based.
Regression clustering methods (sometimes called switching regression methods)
offer an alternative based on probabilistic models. 
Some examples are~\cite{Desarbo88,Quandt78} while the original idea is due 
to~\cite{Bock69}.
However, it has been observed that discarding some dimensions may yield
instabilities in presence of outliers or on small datasets.
For this reason, the method proposed in this paper does not assume that there exist irrelevant dimensions and therefore does not discard any dimensions,
but it models the smallest variances by a single parameter.
Methods based on mixtures of factor analyzers~\cite{McLachlan03,Tipping99}
rely on a latent variables model and on an EM based procedure to cluster
high-dimensional data. More recently, Bocci \emph{et al.}~\cite{Bocci06a}
proposed a similar approach to cluster dissimilarity data.
The model of these methods is a mixture of 
Gaussian densities where the number of parameters is controlled
through the dimension of the latent factor space.
The advantage of such a model is to capture correlations without
estimating full covariance matrices and without dimension
truncation. In this paper, we propose an unified approach for subspace
clustering in the Gaussian mixture model framework which encompasses
these approaches and involves additional regularizations as in parsimonious
models. A precise comparison between our approach and the
mixtures of factor analyzers is achieved
in paragraph~\ref{submodels}.

\section{\label{sec:A-Gaussian-mixture}A Gaussian model for high-dimensional
data}

Clustering divides a given dataset $\{ x_{1},...,x_{n}\}$ of $n$
data points in $\mathbb{R}^{p}$ into $k$ homogeneous groups (see~\cite{Jain99}
for a review). A popular clustering technique uses Gaussian mixture
models, which assume that each class is represented by a Gaussian
probability density. Data are therefore modeled by a density: \begin{equation}
f(x,\theta)=\sum_{i=1}^{k}\pi_{i}\phi(x,\theta_{i}),\label{eq:mixture-model}\end{equation}
 where $\phi$ is a $p$-variate normal density with parameter $\theta_{i}=\{\mu_{i},\Sigma_{i}\}$
and $\pi_{i}$ are the mixing proportions. This model requires to
estimate full covariance matrices and therefore the number of parameters
increases with the square of the dimension. However, due to the \emph{empty
space} phenomenon~\cite{Scott83} it can be assumed that high-dimensional
data live around subspaces with a dimension lower than the one of
the original space. We therefore introduce low-dimensional
class-specific subspaces in order to limit the number of parameters to
estimate.

\subsection{\label{sub:The-new-model}The Gaussian model
$[a_{ij}b_{i}Q_{i}d_{i}]$}

\begin{figure}
\begin{center}\psfrag{Ei}{\footnotesize$E_i$}

\psfrag{xx}{\footnotesize$x$}

\psfrag{Pi}{\footnotesize$P_i(x)$}

\psfrag{mu}{\footnotesize$\mu_i$}

\psfrag{d_1}{\footnotesize$d(x,E_i)$}

\psfrag{d2}{}

\psfrag{d3}{\footnotesize$d(\mu_i,P_i(x))$} 

\psfrag{Ei_perp}{\footnotesize${E_i}^\perp$} 

\psfrag{Pi_perp}{\footnotesize\hspace{-2ex}${P_i}^\perp(x)$} 

\includegraphics[scale=0.45]{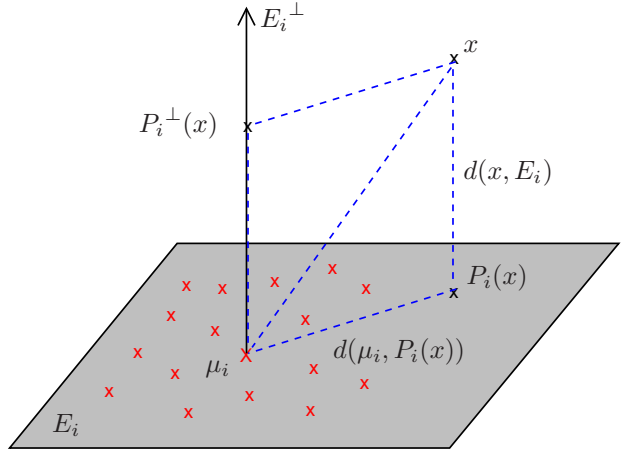}\end{center}

\caption{\label{cap:Illustration_model}The subspaces $\mathbb{E}_{i}$ and
$\mathbb{E}_{i}^{\perp}$ of the $i$th mixture component.}
\end{figure}

As in the classical Gaussian mixture model framework, we assume that
class conditional densities are Gaussian $\mathcal{N}_p(\mu_{i},\Sigma_{i})$
with means $\mu_{i}$ and covariance matrices $\Sigma_{i}$, for $i=1,...,k$.
Let $Q_{i}$ be the orthogonal matrix with the eigenvectors of
$\Sigma_{i}$ as columns.
The class conditional covariance matrix $\Delta_{i}$ is therefore
defined in the eigenspace of $\Sigma_{i}$ by: \begin{equation}
\Delta_{i}=Q_{i}^{t}\,\Sigma_{i}\, Q_{i}.\label{eq:Definition-Delta-i}\end{equation}
 The matrix $\Delta_{i}$ is thus a diagonal matrix which contains
the eigenvalues of $\Sigma_{i}$. It is further assumed that $\Delta_{i}$
is divided into two blocks: 
\begin{equation}
\label{eqdelta}
\small \Delta_i=\left(  \begin{array}{c@{}c} \begin{array}{|ccc|}\hline a_{i1} & & 0\\ & \ddots & \\ 0 & & a_{id_i}\\ \hline \end{array} & \mathbf{0}\\ \mathbf{0} &  \begin{array}{|cccc|}\hline b_i & & & 0\\ & \ddots & &\\  & & \ddots &\\ 0 & & & b_i\\ \hline \end{array} \end{array}\right)  \begin{array}{cc} \left.\begin{array}{c} \\\\\\\end{array}\right\}  & d_{i}\vspace{1.5ex}\\ \left.\begin{array}{c} \\\\\\\\\end{array}\right\}  & (p-d_{i})\end{array}
\end{equation}
with $a_{ij}>b_{i}$, $j=1,...,d_{i}$, and where $d_i\in
\{1,\dots,p-1\}$ is unknown.
The class specific subspace $\mathbb{E}_{i}$ 
is defined as the affine space spanned by the $d_{i}$
eigenvectors associated to the eigenvalues $a_{ij}$ and such that
$\mu_{i}\in\mathbb{E}_{i}$.
Similarly, the affine subspace $\mathbb{E}_{i}^{\perp}$
is such that $\mathbb{E}_{i}\oplus\mathbb{E}_{i}^{\perp}=\mathbb{R}^{p}$
and $\mu_i\in\mathbb{E}_{i}^{\perp}$. 
In this subspace $\mathbb{E}_{i}^{\perp}$, the variance is modeled
by the single parameter $b_{i}$. Let $P_{i}(x)=\tilde{Q_{i}}\tilde{Q_{i}}^{t}(x-\mu_{i})+\mu_{i}$
and $P_{i}^{\perp}(x)=\bar{Q}_i\bar{Q}_i^t(x-\mu_{i})+\mu_{i}$
be the projection of $x$ on $\mathbb{E}_{i}$ and $\mathbb{E}_{i}^{\perp}$
respectively, where $\tilde{Q_{i}}$ is made of the $d_{i}$ first
columns of $Q_{i}$ supplemented by $(p-d_i)$
zero columns and $\bar{Q}_i=(Q_i-\tilde{Q}_i)$.
Thus, $\mathbb{E}_{i}$ is called the specific subspace of the $i$th group since most of the data live on or near this subspace.
In addition, the dimension $d_i$ of the subspace $\mathbb{E}_{i}$ can be considered as the instrinsic dimension of the $i$th group, \emph{i.e.} the number of dimensions required to describe the main features of this group.
Figure~\ref{cap:Illustration_model} summarizes these notations. Following the
notation system of~\cite{Celeux95}, our mixture model is 
denoted by $[a_{ij}b_{i}Q_{i}d_{i}]$ in the sequel.

\subsection{The sub-models of $[a_{ij}b_{i}Q_{i}d_{i}]$}
\label{submodels}

By fixing some parameters to be common within or between classes,
we obtain particular models which correspond to different regularizations.
In the following, ``free $Q_{i}$''
means that $Q_{i}$ is specific for each class $C_{i}$ and {}``common
$Q_{i}$'' means that for each $i=1,...,k$, $Q_{i}=Q$ and consequently
the class orientations are the same. The family~$[a_{ij}b_{i}Q_{i}d_{i}]$
is divided into three categories: models with free orientations, 
common orientations and common covariance matrices.

\paragraph{Models with free orientations}

They assume that the groups live in subspaces
with different orientations, \emph{i.e.} the matrices $Q_{i}$ are
specific to each group. Clearly, the general model $[a_{ij}b_{i}Q_{i}d_{i}]$
belongs to this category. Fixing the dimensions $d_{i}$ to be common between
the classes yields the model~$[a_{ij}b_{i}Q_{i}d]$ which corresponds
to the model of~\cite{Tipping99}.
Indeed, the covariance model given by~(\ref{eq:Definition-Delta-i})
and~(\ref{eqdelta}) can be rewritten as
$\Sigma_i=B_i B_i^t+ D_i$ with $D_i=b_i I_p$, $B_i=Q_i T_i$
and where we have defined 
$$
\small T_i=\left(  \begin{array}{c@{}c} \begin{array}{|ccc|}\hline \sqrt{a_{i1}-b_i} & & 0\\ & \ddots & \\ 0 & & \sqrt{a_{id_i}-b_i}\\ \hline \end{array} & \\\\ \mathbf{0} \\&  \end{array}\right) 
 \begin{array}{cc} \left.\begin{array}{c} \\\\\\\end{array}\right\}  &
d_{i}\vspace{1.5ex}\\ \left.\begin{array}{c} \\\\\\\end{array}\right\}  &
(p-d_{i})\end{array}.
$$
As a consequence, our approach encompasses the mixtures of probabilistic
principal component analysis introduced in~\cite{Tipping99} and
extended in~\cite{McLachlan03} 
to more general matrices $D_i$. In our model, 
$d_i$, the number of columns of $T_i$, depends on the class.
This permits the modeling of a dependence between the number of
factors and the class. Moreover, as illustrated in paragraph~\ref{submodels},
our approach can be combined with a ``parsimonious models'' strategy
to further limit the number of parameters to estimate.
Up to our knowledge, this has not been achieved yet in 
the mixture of factor analyzers model.
For instance, if we further assume that $d_{i}=(p-1)$
for all $i=1,...,k$, the model $[a_{ij}b_{i}Q_{i}d_{i}]$ reduces
to the classical GMM with full covariance matrices for each mixture
component which yields in the supervised framework the well known
Quadratic Discriminant Analysis. It is possible to add constraints
on the different parameters to obtain more regularized models. 
Fixing the first $d_{i}$ eigenvalues to be common within
each class, we obtain the more restricted model $[a_{i}b_{i}Q_{i}d_{i}]$.
The model $[a_{i}b_{i}Q_{i}d_{i}]$ often gives satisfying results,
\emph{i.e.} the assumption that each matrix $\Delta_{i}$ contains
only two different eigenvalues, $a_{i}$ and $b_{i}$, seems to be
an efficient way to regularize the estimation of $\Delta_{i}$. Another
type of regularization is to fix the parameters $b_{i}$ to be common
between the classes. This yields the model $[a_{i}bQ_{i}d_{i}]$
which assumes that the variance outside the class-specific subspaces
is common. This can be viewed as modeling the noise in~$\mathbb{E}_{i}^{\perp}$
by a single parameter $b$ which is natural when the data are
obtained in a common acquisition process. This category of models
contains also the models $[ab_{i}Q_{i}d_{i}]$, $[abQ_{i}d_{i}]$
and all models with free $Q_{i}$ and common $d_{i}$.

\paragraph{Models with common orientations}

It is also possible to assume that the class orientations are common,
\emph{i.e.} $Q_{i}=Q$ for each $i=1,...,k$. However, this assumption
does not necessarily imply that the class-specific subspaces are the
same. Indeed, if the dimensions $d_{i}$ are free, the intersection
of the $k$ class-specific subspaces is the one of the class with
the smallest intrinsic dimension. This assumption can be interesting
to model groups with some common properties and with additional specific
characteristics. Several models of this category require
a complex iterative estimation based on the FG algorithm~\cite{Flury86}
and therefore they will be not considered here. Consequently, only
the models $[a_{i}b_{i}Qd]$, $[ab_{i}Qd]$ and $[a_{i}bQd]$ will
be considered in this paper since their parameters can be estimated
using a simple iterative procedure. Note that a model similar to $[a_{ij}bQd]$
was considered by Flury \emph{et al.} in~\cite{Flury97} in the supervised
framework with an additional assumption on the means.

\paragraph{Models with common covariance matrices}

This branch of the family only includes two models $[a_{j}bQd]$
and $[abQd]$. Both models indeed assume that the classes have the
same covariance matrix $\Sigma=Q\Delta Q^{t}$. Particularly, fixing
$d=(p-1)$, the model $[a_{j}bQd]$ reduces to a Gaussian
mixture model (denoted by {}``Com-GMM'' in the following) which
yields in the supervised framework the well known Linear Discriminant
Analysis (LDA). Remark that if $d<(p-1)$, the model $[a_{j}bQd]$
can be viewed as the a combination of a dimension reduction technique
with a GMM with common covariance matrices, but without losing information
since the information carried by the smallest eigenvalues is not discarded.

\subsection{Characteristics of the models}

\begin{table}[p]
\begin{singlespace}
\begin{center}\footnotesize{\begin{tabular}{lcccc}
\hline 
Model&
\begin{tabular}{c}
Number of\tabularnewline
parameters\tabularnewline
\end{tabular}&
\begin{tabular}{c}
Asymptotic\tabularnewline
order\tabularnewline
\end{tabular}&
\begin{tabular}{c}
Nb of prms $k=4$,\tabularnewline
$d=10$, $p=100$\tabularnewline
\end{tabular}&
\begin{tabular}{c}
ML\tabularnewline
estimation\tabularnewline
\end{tabular}\tabularnewline
\hline
$[a_{ij}b_{i}Q_{i}d_{i}]$&
$\rho+\bar{\tau}+2k+D$&
$kpd$&
4231&
CF\tabularnewline
$[a_{ij}bQ_{i}d_{i}]$&
$\rho+\bar{\tau}+k+D+1$&
$kpd$&
4228&
CF\tabularnewline
$[a_{i}b_{i}Q_{i}d_{i}]$&
$\rho+\bar{\tau}+3k$&
$kpd$&
4195&
CF\tabularnewline
$[ab_{i}Q_{i}d_{i}]$&
$\rho+\bar{\tau}+2k+1$&
$kpd$&
4192&
CF\tabularnewline
$[a_{i}bQ_{i}d_{i}]$&
$\rho+\bar{\tau}+2k+1$&
$kpd$&
4192&
CF\tabularnewline
$[abQ_{i}d_{i}]$&
$\rho+\bar{\tau}+k+2$&
$kpd$&
4189&
CF\tabularnewline
$[a_{ij}b_{i}Q_{i}d]$&
$\rho+k(\tau+d+1)+1$&
$kpd$&
4228&
CF\tabularnewline
$[a_{j}b_{i}Q_{i}d]$&
$\rho+k(\tau+1)+d+1$&
$kpd$&
4198&
CF\tabularnewline
$[a_{ij}bQ_{i}d]$&
$\rho+k(\tau+d)+2$&
$kpd$&
4225&
CF\tabularnewline
$[a_{j}bQ_{i}d]$&
$\rho+k\tau+d+2$&
$kpd$&
4195&
CF\tabularnewline
$[a_{i}b_{i}Q_{i}d]$&
$\rho+k(\tau+2)+1$&
$kpd$&
4192&
CF\tabularnewline
$[ab_{i}Q_{i}d]$&
$\rho+k(\tau+1)+2$&
$kpd$&
4189&
CF\tabularnewline
$[a_{i}bQ_{i}d]$&
$\rho+k(\tau+1)+2$&
$kpd$&
4189&
CF\tabularnewline
$[abQ_{i}d]$&
$\rho+k\tau+3$&
$kpd$&
4186&
CF\tabularnewline
\hline 
$[a_{ij}b_{i}Qd_{i}]$&
$\rho+\tau+D+2k$&
$pd$&
1396&
FG\tabularnewline
$[a_{ij}bQd_{i}]$&
$\rho+\tau+D+k+1$&
$pd$&
1393&
FG\tabularnewline
$[a_{i}b_{i}Qd_{i}]$&
$\rho+\tau+3k$&
$pd$&
1360&
FG\tabularnewline
$[a_{i}bQd_{i}]$&
$\rho+\tau+2k+1$&
$pd$&
1357&
FG\tabularnewline
$[ab_{i}Qd_{i}]$&
$\rho+\tau+2k+1$&
$pd$&
1357&
FG\tabularnewline
$[abQd_{i}]$&
$\rho+\tau+k+2$&
$pd$&
1354&
FG\tabularnewline
$[a_{ij}b_{i}Qd]$&
$\rho+\tau+kd+k+1$&
$pd$&
1393&
FG\tabularnewline
$[a_{j}b_{i}Qd]$&
$\rho+\tau+k+d+1$&
$pd$&
1363&
FG\tabularnewline
$[a_{ij}bQd]$&
$\rho+\tau+kd+2$&
$pd$&
1390&
FG\tabularnewline
$[a_{i}b_{i}Qd]$&
$\rho+\tau+2k+1$&
$pd$&
1357&
IP\tabularnewline
$[ab_{i}Qd]$&
$\rho+\tau+k+2$&
$pd$&
1354&
IP\tabularnewline
$[a_{i}bQd]$&
$\rho+\tau+k+2$&
$pd$&
1354&
IP\tabularnewline
\hline 
$[a_{j}bQd]$&
$\rho+\tau+d+2$&
$pd$&
1360&
CF\tabularnewline
$[abQd]$&
$\rho+\tau+3$&
$pd$&
1351&
CF\tabularnewline
\hline 
Full-GMM&
$\rho+kp(p+1)/2$&
$kp^{2}/2$&
20603&
CF\tabularnewline
Com-GMM&
$\rho+p(p+1$)/2&
$p^{2}/2$&
5453&
CF\tabularnewline
Diag-GMM&
$\rho+kp$&
$2kp$&
803&
CF\tabularnewline
Sphe-GMM&
$\rho+k$&
$kp$&
407&
CF\tabularnewline
\hline
\end{tabular}}\end{center}
\end{singlespace}

\caption{\label{cap:Some-properties-of-HDDA-models}Properties of the HDDC
models: $\rho=kp+k-1$ is the number of parameters required for the
estimation of means and proportions, $\bar{\tau}=\sum_{i=1}^{k}d_{i}[p-(d_{i}+1)/2]$
and $\tau=d[p-(d+1)/2]$ are the number of parameters required for
the estimation of $\tilde{Q_{i}}$ and $\tilde{Q}$, and $D=\sum_{i=1}^{k}d_{i}$.
For asymptotic orders, we assume that $k\ll d\ll p$. CF means that
the ML estimates are closed form. IP means that the ML estimation
needs an iterative procedure. FG means that the ML estimation requires
the iterative FG algorithm.}
\end{table}

Our family of models presented above only requires the estimation
of $d_{i}$-dimensional subspaces and therefore the different models
are significantly more parsimonious than the general Gaussian model
if $d_{i}\ll p$. Table~\ref{cap:Some-properties-of-HDDA-models}
summarizes some properties of the models considered here. The second
column of this table gives the number of parameters to estimate. The
third column provides the asymptotic order of the number of parameters
(\emph{i.e.} with the assumption that $k\ll d\ll p$). The fourth
column gives the number of parameters for the particular case $k=4$,
$p=100$ and $\forall i,\, d_{i}=10$. The last column indicates whether
the Maximum Likelihood~(ML) updates are in closed form or not. These
characteristics are also given for five Gaussian mixture models: GMM
with full covariance matrices for each class (Full-GMM), with common
covariance matrices between classes (Com-GMM), with diagonal covariance
matrices (Diag-GMM), with spherical covariance matrices (Sphe-GMM).
Note that Celeux and Govaert recommend in~\cite{Celeux95} to make
use of the models Diag-GMM and Sphe-GMM in clustering problems. We
can observe that all models of our family require the estimation of
fewer parameters than both Full-GMM and Com-GMM. In the particular
case of 100-dimensional data, made of 4 classes and with common intrinsic
dimensions $d_{i}$ equal to 10, the model $[a_{ij}b_{i}Q_{i}d_{i}]$
only requires the estimation of 4~231 parameters whereas Full-GMM
and Com-GMM requires respectively the estimation of 20~603 and 5~453
parameters. Remark that the model $[a_{ij}b_{i}Q_{i}d_{i}]$, which
gives rise to quadratic separation between the groups, requires the
estimation of fewer parameters than Com-GMM, which gives rise to linear
separation between the groups.

\section{\label{sec:Clustering-method}High-dimensional data clustering}

In this section, we derive the EM-based clustering framework for the
model $[a_{ij}b_{i}Q_{i}d_{i}]$ and its sub-models. The related
clustering method is denoted
by High-Dimensional Data Clustering (HDDC). Let us recall that unsupervised
classification organizes data in homogeneous groups using only the
observed values of the $p$ explanatory variables. Usually, in model-based
clustering, the parameters $\theta=\{\pi_{1},...,\pi_{k},\theta_{1},...,\theta_{k}\}$
with $\theta_{i}=\{\mu_{i},\Sigma_{i}\}$ are estimated by the EM
algorithm which repeats iteratively E and M steps. The reader could
refer to~\cite{Mclachlan97} for further informations on the EM algorithm
and its extensions. In particular, the models presented in this paper
can be also used in the Classification EM and Stochastic EM
algorithms~\cite{Celeux92}. Using our parameterization,
the EM algorithm for estimating
$\theta=\{\pi_{i},\mu_{i},\Sigma_{i},a_{ij},b_{i},Q_{i},d_{i}\}$
is detailed in the following.

\subsection{\label{sub:The-E-step}The E step}

This step computes, at iteration $q$ and for each $i=1,...,k$ and
$j=1,...,n$, the conditional probability $t_{ij}^{(q)}=\mathbb{P}(x_{j}\in C_{i}^{(q-1)}|x_{j})$
which can be written from~(\ref{eq:mixture-model}) and using the
Bayes formula as follows:\[
t_{ij}^{(q)}={\pi_{i}^{(q-1)}\phi(x_{j},\theta_{i}^{(q-1)})}\left/{\sum_{\ell=1}^{k}\pi_{\ell}^{(q-1)}\phi(x_{j},\theta_{\ell}^{(q-1)})}\right..\]
Note that this conditional probability 
is mainly
based on $\pi_{i}^{(q-1)}\phi(x_{j},\theta_{i}^{(q-1)})$.
and thus can be rewritten using the
parameters of the model $[a_{ij}b_{i}Q_{i}d_{i}]$. 
In order not to overload the equations, 
the index of the current iteration $q$ is omitted in the remainder of this
paragraph. Writing $\phi(x,\theta_{i})$ with the new class conditional
covariance matrix $\Delta_{i}$, we obtain:\[
-2\log(\phi(x,\theta_{i}))=(x-\mu_{i})^{t}(Q_{i}\Delta_{i}Q_{i}^{t})^{-1}(x-\mu_{i})+\log(\det\Delta_{i})+p\log(2\pi).\]
Since $Q_i^tQ_i=I_p$ and $Q_i=\tilde{Q}_i+\bar{Q}_i$, the above 
matrix inverse can be expanded as $(Q_{i}\Delta_{i}Q_{i}^{t})^{-1}= \tilde{Q}_i\Delta_i^{-1}\tilde{Q}_i^t + \bar{Q}_i\Delta_i^{-1}\bar{Q}_i^t $ and thus:\begin{eqnarray*}
-2\log(\phi(x,\theta_{i})) & = & (x-\mu_{i})^{t}\tilde{Q}_{i}\Delta_{i}^{-1}\tilde{Q}_{i}^{t}(x-\mu_{i})+(x-\mu_{i})^{t}\bar{Q_{i}}\Delta_{i}^{-1}\bar{Q}_{i}^{t}(x-\mu_{i})\\
 &  & +\log(\det\Delta_{i})+p\log(2\pi).\end{eqnarray*}
Taking into account the structure of $\Delta_{i}$ and using the relations $\tilde{Q}_i (\tilde{Q}_i^t\tilde{Q}_i) = \tilde{Q}_i$ and $\bar{Q}_i \left(\bar{Q}_i^t\bar{Q}_i\right) = \bar{Q}_i $, it yields:\[
-2\log(\phi(x,\theta_{i}))=\|\tilde{Q}_{i}\tilde{Q}_{i}^t(x-\mu_{i})\|_{\mathcal{A}_{i}}^
{2}+\frac{1}{b_{i}}\|\bar{Q}_{i}\bar{Q}_{i}^t(x-\mu_{i})\|^{2}
+\log(\det\Delta_{i})+p\log(2\pi),\]
where $\|.\|_{\mathcal{A}_{i}}^{2}$  is the norm on $\mathbb{E}_{i}$ such
as $\| x\|_{\mathcal{A}_{i}}^{2}=x^{t}\mathcal{A}_{i}x$ with
$\mathcal{A}_{i}=\tilde{Q_{i}}\Delta_{i}^{-1}\tilde{Q_{i}}^{t}$.
From the definitions of $P_{i}$ and $P_{i}^{\perp}$ (Paragraph~\ref{sub:The-new-model}) and in view of Figure~\ref{cap:Illustration_model},
we have:\[
-2\log(\phi(x,\theta_{i}))=\|\mu_{i}-P_{i}(x)\|_{\mathcal{A}_{i}}^{2}+\frac{1}{
b_{i}}\| x-P_{i}(x)\|^{2}+\log(\det\Delta_{i})+p\log(2\pi).\]
The relation
$\log(\det\Delta_{i})=\sum_{j=1}^{d_{i}}\log(a_{ij})+(p-d_{i})\log(b_{i})$
allows to conclude that:\[
t_{ij}=1\left/\sum_{\ell=1}^{k}\exp\left(\frac{1}{2}(K_{i}(x_{j})-K_{\ell}(x_{j}
))\right)\right.,\]
 where $K_{i}(x)=-2\log(\pi_{i}\phi(x,\theta_{i}))$ is called the
cost function and is defined by:\[
K_{i}(x)=\Vert\mu_{i}-P_{i}(x)\|_{\mathcal{A}_{i}}^{2}+\frac{1}{b_{i}}\Vert
x-P_{i}(x)\Vert^{2}+\sum_{j=1}^{d_{i}}\log(a_{ij})+(p-d_{i})\log(b_{i}
)-2\log(\pi_{i}).\label{eq:definition-of-Ki}\]
Let us note that $K_{i}(x)$ is mainly based on two distances:
the distance between the projection of $x$ on $\mathbb{E}_{i}$ and
the mean of the class and the distance between the observation and
the subspace $\mathbb{E}_{i}$. This cost function favors the assignment
of a new observation to the class for which it is close to the subspace
and for which its projection on the class subspace is close to the
mean of the class. The variance terms $a_{ij}$ and $b_{i}$ balance
the importance of both distances. For example, if the data are very
noisy, \emph{i.e.} $b_{i}$ is large, it is natural to balance the
distance $\Vert x-P_{i}(x)\Vert^{2}$ by $1/b_{i}$ in order to take
into account the large variance in $\mathbb{E}_{i}^{\perp}$.

\subsection{\noindent \label{sub:The-M-step}The M step}

This step maximizes at iteration $q$ the conditional likelihood and
uses the following update formulas. Mixture proportions and means
are estimated by:\[
\hat{\pi}_{i}^{(q)}=\frac{n_{i}^{(q)}}{n},\,\,\hat{\mu}_{i}^{(q)}=\frac{1}{n_{i}
^{(q)}}\sum_{j=1}^{n}t_{ij}^{(q)}x_{j},\]
 where $n_{i}^{(q)}=\sum_{j=1}^{n}t_{ij}^{(q)}$. Moreover, the update
formula for the empirical covariance matrix of the fuzzy class $C_{i}$
is:\[
W_{i}^{(q)}=\frac{1}{n_{i}^{(q)}}\sum_{j=1}^{n}t_{ij}^{(q)}(x_{j}-\hat{\mu}_{i}^
{(q)})(x_{j}-\hat{\mu}_{i}^{(q)})^{t}.\]
 The estimation of the specific parameters of HDDC is detailed below.
Proofs of the following results are given in the Appendix.

\paragraph{Models with free orientations}

The ML estimators of model parameters are closed form for this category
of models.

\noindent -- Subspace $\mathbb{E}_{i}$: the $d_{i}$ first columns
of $Q_{i}$ are estimated by the eigenvectors associated with the
$d_{i}$ largest eigenvalues $\lambda_{ij}$ of $W_{i}$.

\noindent -- Model $[a_{ij}b_{i}Q_{i}d_{i}]$: the estimator of $a_{ij}$
is $\hat{a}_{ij}=\lambda_{ij}$ and the estimator of $b_{i}$ is the
mean of the $(p-d_{i})$ smallest eigenvalues of $W_{i}$ and can
be written as follows:\begin{equation}
\hat{b}_{i}=\frac{1}{(p-d_{i})}\left(\trace(W_{i})-\sum_{j=1}^{d_{i}}\lambda_{ij
}\right).\label{eq:estimator_of_b_i}\end{equation}
-- Model $[a_{ij}bQ_{i}d_{i}]$: the estimator of $a_{ij}$ is
$\hat{a}_{ij}=\lambda_{ij}$
and the estimator of $b$ is:\begin{equation}
\hat{b}=\frac{1}{(p-\xi)}\left(\trace(W)-\sum_{i=1}^{k}\hat{\pi}_{i}\sum_{j=1}^{
d_{i}}\lambda_{ij}\right),\label{eq:estimator_of_b}\end{equation}
where $\xi=\sum_{i=1}^{k}\hat{\pi}_{i}d_{i}$ and
$W=\sum_{i=1}^{k}\hat{\pi}_{i}W_{i}$
is the within-covariance matrix.

\noindent -- Model $[a_{i}b_{i}Q_{i}d_{i}]$: the estimator of $b_{i}$
is given by (\ref{eq:estimator_of_b_i}) and the estimator of $a_{i}$
is:\begin{equation}
\hat{a}_{i}=\frac{1}{d_{i}}\sum_{j=1}^{d_{i}}\lambda_{ij}.
\label{eq:estimator_of_a_i}\end{equation}
 -- Model $[ab_{i}Q_{i}d_{i}]$: the estimator of $b_{i}$ is given
by (\ref{eq:estimator_of_b_i}) and the estimator of $a$ is:\begin{equation}
\hat{a}=\frac{1}{\xi}\sum_{i=1}^{k}\hat{\pi}_{i}\sum_{j=1}^{d_{i}}\lambda_{ij}
.\label{eq:estimator_of_a}\end{equation}
 -- Model $[a_{i}bQ_{i}d_{i}]$: estimators of $a_{i}$ and $b$
are respectively given by (\ref{eq:estimator_of_a_i}) and
(\ref{eq:estimator_of_b}).\\
-- Model $[abQ_{i}d_{i}]$: estimators of $a$ and $b$ are respectively
given by (\ref{eq:estimator_of_a}) and (\ref{eq:estimator_of_b}).

\noindent -- Models with common dimensions: the estimators of the
models with common dimensions $d_{i}$ can be obtained from the previous
ones by replacing the values $d_{i}$ by $d$ for each $i=1,...,k$.
In this case, equations~(\ref{eq:estimator_of_b}) and~(\ref{eq:estimator_of_a})
can be simplified as:\begin{gather}
\hat{a}=\frac{1}{d}\sum_{j=1}^{d}\lambda_{j},\label{eq:estimator_of_a_bis}\\
\hat{b}=\frac{1}{(p-d)}\left(\trace(W)-\sum_{j=1}^{d}\lambda_{j}\right),
\label{eq:estimator_of_b_bis}\end{gather}
where $\lambda_{j}$ is the $j$th largest eigenvalue of $W$.

\noindent -- Model $[a_{j}b_{i}Q_{i}d]$: the estimator of $a_{j}$
is $\hat{a}_{j}=\lambda_{j}$ and the estimator of $b_{i}$ is
(\ref{eq:estimator_of_b_i}).

\noindent -- Model $[a_{j}bQ_{i}d]$: the estimator of $a_{j}$ is
$\hat{a}_{j}=\lambda_{j}$ and the estimator of $b$ is (\ref{eq:estimator_of_b_bis}).

\paragraph{Models with common orientations}

Here, we assume in addition that the dimensions $d_{i}$ are common
between classes. The following ML estimators require an iterative
procedure.

\noindent -- Subspace $\mathbb{E}_{i}$: Given $a_{i}$ and $b_{i}$,
the $d$ first columns of $Q$ are estimated by the eigenvectors associated
to the $d$ largest eigenvalues of the matrix $M$ defined by:\[
M(a_{1},...,a_{k},b_{1},...,b_{k})=\sum_{i=1}^{k}n_{i}(\frac{1}{b_{i}}-\frac{1}{
a_{i}})W_{i}.\]
\noindent -- Model $[a_{i}b_{i}Qd]$: given $Q$, estimators of
$a_{i}$ and $b_{i}$ are:\begin{gather}
\hat{a}_{i}(Q)=\frac{1}{d}\sum_{j=1}^{d}q_{j}^{t}W_{i}q_{j},
\label{eq:estimator_of_a_iQ}\\
\hat{b}_{i}(Q)=\frac{1}{(p-d)}\left(\trace(W_{i})-\sum_{j=1}^{d}q_{j}^{t}W_{i}q_
{j}\right).\label{eq:estimator_of_b_iQ}\end{gather}
-- Model $[a_{i}bQ_{i}d_{i}]$: given $Q$, the estimator of $a_{i}$
is (\ref{eq:estimator_of_a_iQ}) and the estimator of $b$
is:\begin{equation}
\hat{b}(Q)=\frac{1}{(p-d)}\left(\trace(W)-\sum_{j=1}^{d}q_{j}^{t}Wq_{j}
\right).\label{eq:estimator_of_bQ}\end{equation}
\noindent -- Model $[ab_{i}Qd]$: given $Q$, the estimator of $b_{i}$
is (\ref{eq:estimator_of_b_iQ}) and the estimator of $a$
is:\begin{equation}
\hat{a}(Q)=\frac{1}{d}\sum_{j=1}^{d}q_{j}^{t}Wq_{j}.\label{eq:estimator_of_aQ}
\end{equation}
\noindent -- Model $[a_{i}bQd]$: given $Q$, estimators of $a_{i}$
and $b$ are respectively (\ref{eq:estimator_of_a_iQ}) and
(\ref{eq:estimator_of_bQ}).

For example, it is possible to use the following iterative procedure
to estimate the parameters associated to the model $[a_{i}b_{i}Qd]$: 

\noindent -- Initialization: the $d$ first columns of $Q^{(0)}$
are the eigenvectors associated with the $d$ largest eigenvalues
of $W$. 

\noindent -- Until convergence: $a_{i}^{(\ell)}=\hat{a}_{i}(Q^{(\ell-1)})$,
$b_{i}^{(\ell)}=\hat{b}_{i}(Q^{(\ell-1)})$ and the $d$ first columns
of $Q^{(\ell)}$ are the eigenvectors associated to the $d$ largest
eigenvalues of the matrix
$M(a_{1}^{(\ell)},...,a_{k}^{(\ell)},b_{1}^{(\ell)},...,b_{k}^{(\ell)})$.

\paragraph{Models with common covariance matrices}

In this category of models, the parameters can be updated in closed
form.

\noindent -- Subspace $\mathbb{E}_{i}$: the $d$ first columns of
the matrix $Q$ are the eigenvectors associated to the $d$ largest
eigenvalues of~$W$.

\noindent -- Model $[a_{j}bQd]$: the estimator of $a_{j}$ is
$\hat{a}_{j}=\lambda_{j}$
and the estimator of $b$ is~(\ref{eq:estimator_of_b_bis}).

\noindent -- Model $[abQd]$: estimators of $a$ and $b$ are respectively
given by (\ref{eq:estimator_of_a_bis}) and (\ref{eq:estimator_of_b_bis}).

\subsection{\label{sub:Intrinsic-dimension-estimation}Hyper-parameters estimation}

\begin{figure}
\begin{center}\includegraphics[%
  width=0.49\columnwidth]{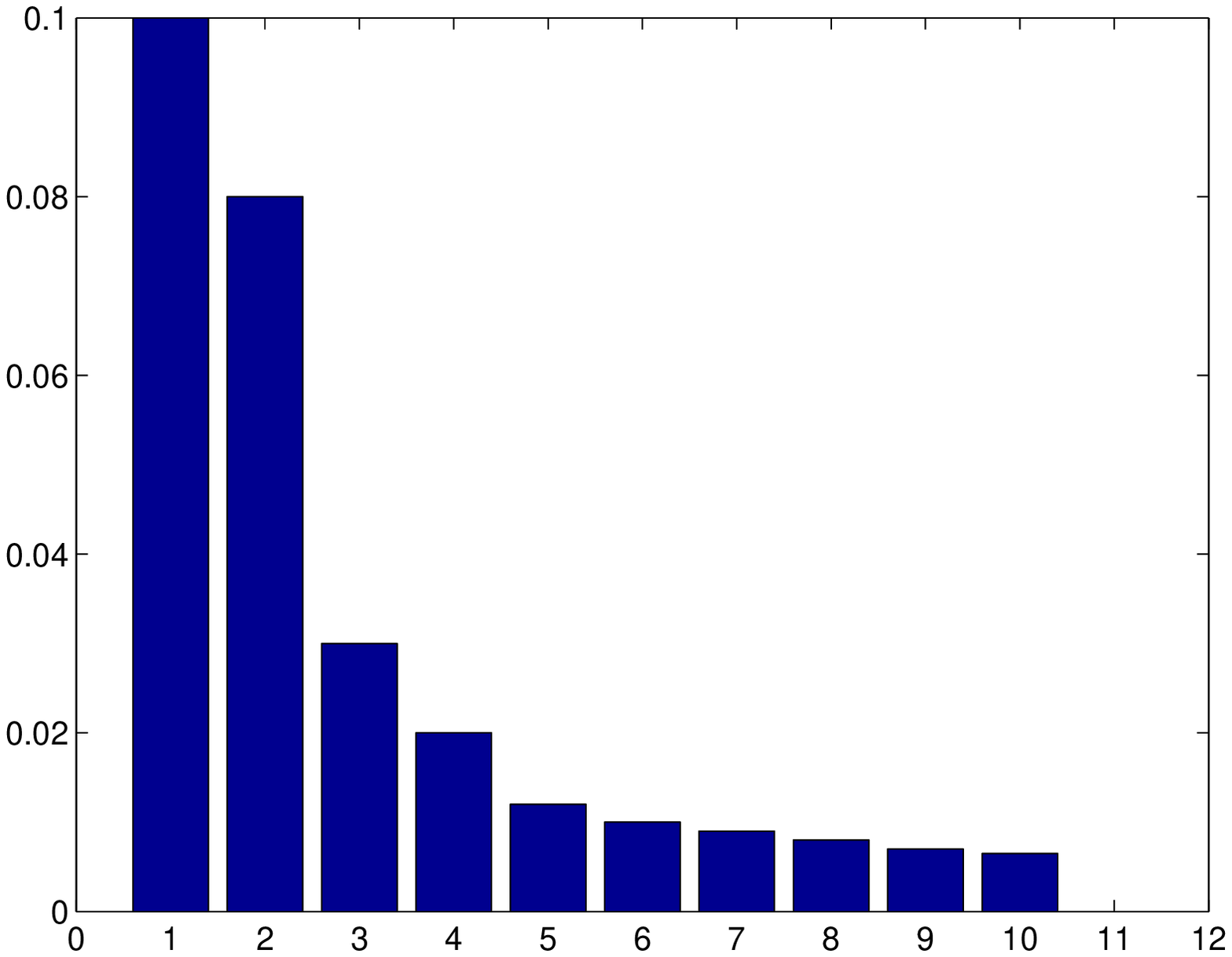}\includegraphics[%
  width=0.49\columnwidth]{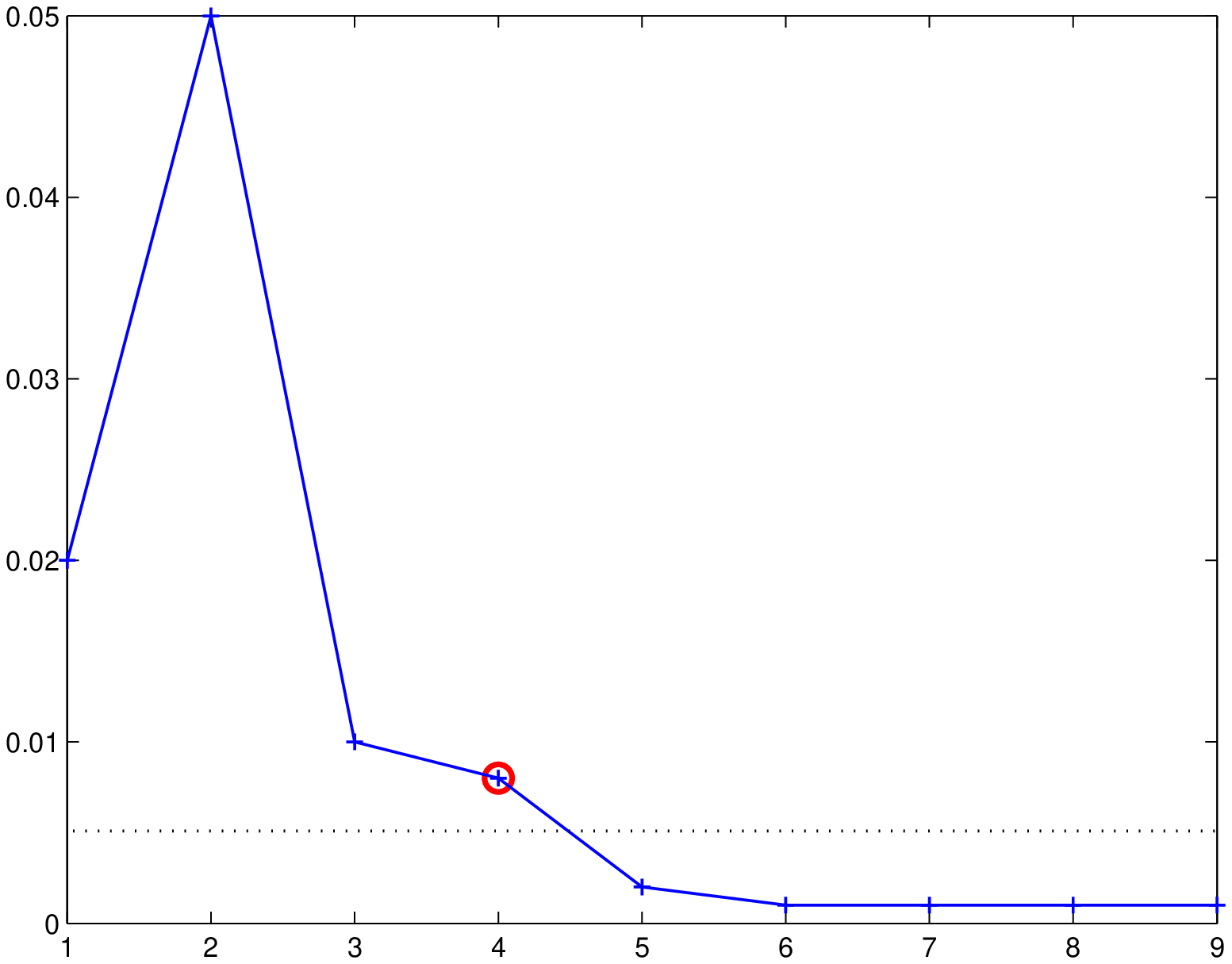}\vspace{-3ex}\end{center}

\caption{\label{cap:Scree-test}Estimation of the intrinsic dimension $d_{i}$
using the scree-test of Cattell: plot of ordered eigenvalues of $\Sigma_{i}$
(left) and differences between consecutive eigenvalues (right).}
\end{figure}

Within the M step, the intrinsic dimensions
of each subclass have to be estimated.
 This is a difficult problem with no unique technique
to use. Our approach is based on the eigenvalues of the class conditional
covariance matrix $\Sigma_{i}$ of the class $C_{i}$. The $j$th
eigenvalue of $\Sigma_{i}$ corresponds to the fraction of the full
variance carried by the $j$th eigenvector of $\Sigma_{i}$.
The class specific dimension $d_{i},\, i=1,...,k$ is estimated through
the scree-test of Cattell~\cite{Cattell66} which looks for a
break in the eigenvalues scree.
The selected dimension is the one for which the subsequent eigenvalues
differences are smaller than a threshold. 
Figure~\ref{cap:Scree-test} illustrates
this method: the graph on the right shows that the differences between
eigenvalues after the fourth one are smaller than the threshold (dashed
line). Thus, in this case, four dimensions will be chosen and this
corresponds indeed to a break in the scree (left graph). In our experiments,
the threshold is chosen using the probabilistic criterion BIC~\cite{Schwarz78}
which consists in minimizing
$\bic(m)=-2\log(L)+\nu(m)\log(n)$,
where $\nu(m)$ is the number of parameters of the model $m$ given
in Table~\ref{cap:Some-properties-of-HDDA-models} for HDDC, $L$
is the likelihood and $n$ is the number of observations.
In addition, this approach allows to estimate $k$
parameters by choosing only the value of the threshold $t$. In the case of
common intrinsic dimensions between the groups, the dimension $d$ is directly
determined using BIC. The second hyper-parameter to estimate in any clustering
method is the number of groups $k$. This parameter is also selected
thanks to the BIC criterion, see the experiments presented in
Section~\ref{sec:Experimental-results}.

\subsection{Numerical considerations}

First, it is important to remark that the parametrization of the
Gaussian model proposed here provides an explicit expression of
$\Sigma_{i}^{-1}$ whereas other classical methods, like Full-GMM and Com-GMM,
need to numerically invert empirical covariance matrices which usually fails for singularity reasons.
Some solutions however exist to overcome this problem 
for the models Full-GMM and Com-GMM, see for instance~\cite{Krzanowski95}.
In contrast, this
problem does not arise with HDDC since the cost function $K_{i}$ does not
require to invert $\Sigma_{i}$. 
Moreover, it appears in~(\ref{eq:definition-of-Ki}) that the
cost function $K_{i}$ does not use the projection on the subspace
$\mathbb{E}_{i}^{\perp}$ and consequently does not require the computation of
the $(p-d_i)$ latest columns of the orientation matrix $Q_i$. 
In Section~\ref{sub:The-M-step}, it is shown that the ML estimators of these columns are the eigenvectors associated to the $(p-d_{i})$ smallest eigenvalues of the empirical covariance matrix $W_{i}$. Therefore, HDDC does not depend on these eigenvectors whose determination is numerically unstable. 
Thus, HDDC is robust with respect to ill-conditioning 
and singularity problems.
In addition, it is also
possible to use this feature to reduce computing time by using the
Arnoldi method~\cite{Lehoucq98} which only provides the largest
eigenvalues and the associated eigenvectors of an ill-conditioned matrix. During
our experiments, we noticed a reduction by a factor 60 of the computing time on
a 1024-dimensional dataset compared to the classical approach. Furthermore, in
the special case where the number of observations of a group $n_i$ is smaller
than the dimension $p$, our parametrization allows to use a linear algebra
trick. Indeed, in this case, it is better from a numerical point of view to
compute the eigenvectors of the $n_i \times n_i$ matrix
$\Upsilon_{i}\Upsilon_{i}^{t}$ than those of the $p \times p$ matrix
$\Upsilon_{i}^{t}\Upsilon_{i}$, where $\Upsilon_{i}$ is the $n_i \times p$
matrix containing the mean-centered observations. In the case of
data containing 13 observations in a 1024-dimensional space, it has been
noticed a
reduction by a factor 500 of the computing time compared to the classical
approach.

\section{\label{sec:Experimental-results}Experimental results}

In this section, we present results for artificial and real datasets
illustrating the main features of HDDC. In the following experiments,
HDDC will be compared to 3 classical Gaussian mixture models: GMM
with full covariance matrices for each class (Full-GMM), with diagonal
covariance matrices (Diag-GMM), with spherical covariance matrices
(Sphe-GMM). A numerical regularization was necessary to invert the
covariance matrices in the clustering method associated to the model
Full-GMM, so that it is able to work with data of dimension larger
than 50.

\subsection{Simulation study: model selection}

\begin{table}
\begin{center}\small{\begin{tabular}{|l|cccccc|}
\hline 
Simulated &
\multicolumn{6}{c|}{HDDC model}\tabularnewline
data model&
 $[a_{ij}b_{i}Q_{i}d_{i}]$&
 $[a_{ij}bQ_{i}d_{i}]$&
 $[a_{i}b_{i}Q_{i}d_{i}]$&
 $[a_{i}bQ_{i}d_{i}]$&
 $[ab_{i}Q_{i}d_{i}]$&
 $[abQ_{i}d_{i}]$\tabularnewline
\hline
$[a_{ij}b_{i}Q_{i}d_{i}]$&
357 &
 373 &
 \textbf{349 }&
 359 &
 \textbf{349 }&
 360\tabularnewline
 $[a_{ij}bQ_{i}d_{i}]$&
 403 &
 404 &
 397 &
 \textbf{396 }&
 397 &
 397\tabularnewline
 $[a_{i}b_{i}Q_{i}d_{i}]$&
 389 &
 419 &
 \textbf{377 }&
 391 &
 \textbf{377 }&
 394\tabularnewline
 $[a_{i}bQ_{i}d_{i}]$&
 438 &
 440 &
 \textbf{419 }&
 \textbf{419 }&
 420 &
 420\tabularnewline
 $[ab_{i}Q_{i}d_{i}]$&
 399 &
 433 &
 \textbf{380 }&
 402 &
 384 &
 403\tabularnewline
 $[abQ_{i}d_{i}]$&
 456 &
 451 &
 428 &
 \textbf{427 }&
 434 &
 433 \tabularnewline
\hline
\end{tabular}}\vspace{-2ex}\end{center}

\caption{\label{cap:Table-BIC-2}BIC value for the HDDC models on different
simulated datasets (the best ones are in bold).}
\end{table}

\begin{table}[t]
\begin{center}\small{\begin{tabular}{|l|cccccc|}
\hline 
Simulated&
\multicolumn{6}{c|}{HDDC model}\tabularnewline
data model&
 $[a_{ij}b_{i}Q_{i}d_{i}]$&
 $[a_{ij}bQ_{i}d_{i}]$&
 $[a_{i}b_{i}Q_{i}d_{i}]$&
 $[a_{i}bQ_{i}d_{i}]$&
 $[ab_{i}Q_{i}d_{i}]$&
 $[abQ_{i}d_{i}]$\tabularnewline
\hline
$[a_{ij}b_{i}Q_{i}d_{i}]$&
0.967 &
 0.828 &
 0.973 &
 0.919 &
 \textbf{0.975 }&
 0.903\tabularnewline
 $[a_{ij}bQ_{i}d_{i}]$&
 0.730 &
 0.727 &
 0.779 &
 \textbf{0.782 }&
 0.758 &
 0.751\tabularnewline
 $[a_{i}b_{i}Q_{i}d_{i}]$&
 0.979 &
 0.871 &
 0.983 &
 0.929 &
 \textbf{0.986 }&
 0.917\tabularnewline
 $[a_{i}bQ_{i}d_{i}]$&
 0.826 &
 0.800&
 \textbf{0.882 }&
 0.863 &
 0.875 &
 0.865\tabularnewline
 $[ab_{i}Q_{i}d_{i}]$&
 0.965 &
 0.825 &
 \textbf{0.980 }&
 0.844 &
 0.952 &
 0.822\tabularnewline
 $[abQ_{i}d_{i}]$&
 0.712 &
 0.752 &
 \textbf{0.797 }&
 0.793 &
 0.711 &
 0.707 \tabularnewline
\hline
\end{tabular}}\vspace{-2ex}\end{center}

\caption{\label{cap:Table-BIC-3}Cluster recognition rate for the HDDC models
on different simulated datasets (the best ones are in bold).}
\end{table}

Given that HDDC is a model-based clustering method, the well-known
criterion BIC can be used for selecting the best adapted model to
the data. Here, we used BIC and the cluster recognition rate to compare
the different models of HDDC. The cluster recognition rate can be
computed since true partitions are known and is defined as the maximum rate
over the correct matchings between the true groups and the found clusters.
It is impossible to report in this section numerical experiments for
all the discussed models. Therefore, we limit ourselves to models
with free orientations since we believe that these models are able
to tackle different situations. We performed extensive simulations
(50 replications for each of the 6 data models) and then used the
6 different models with free orientations in HDDC to cluster the simulated
data. For each dataset, 3 Gaussian densities are simulated in $\mathbb{R}^{100}$
according to one of the 6 models with free orientations, \emph{i.e.}
free matrices $Q_{i}$, and with the following parameters: $\{ d_{1},d_{2},d_{3}\}=\{2,5,10\}$,
$\{\pi_{1},\pi_{2},\pi_{3}\}=\{0.4,0.3,0.3\}$ and close means and
random matrices $Q_{i}$. Each one of the 6 datasets was made of 1000
points. Tables~\ref{cap:Table-BIC-2} and~\ref{cap:Table-BIC-3}
present respectively the BIC value and the cluster recognition rate
on average for the 6 considered HDDC models on the different simulated
datasets. First of all, it appears that BIC and the cluster recognition
rate select in general the same models and this confirm that BIC is
a useful tool in model-based clustering. Unsurprisingly, the models
used to simulate the data obtain small BIC values and satisfying cluster
recognition rates. However, it appears that the model $[a_{i}b_{i}Q_{i}d_{i}]$
is usually selected by BIC as the best model and its cluster recognition
rates are very good for each type of simulated data. Thus, the model
$[a_{i}b_{i}Q_{i}d_{i}]$ seems to have the right number of degrees
of freedom and the assumption that $\Delta_{i}$ has only 2 different
eigenvalues is an efficient way to regularize the estimation.
Note that models $[a_{i}bQ_{i}d_{i}]$ and $[ab_{i}Q_{i}d_{i}]$
are also often selected by BIC and provide good cluster recognition
rates.

\subsection{Simulation study: hyper-parameters selection}

\begin{table}
\begin{center}\begin{tabular}{|c|l|c|}
\hline 
Nb of groups $k$&
Dimensions $d_{i}$&
BIC value\tabularnewline
\hline
2&
2,16&
414\tabularnewline
\textbf{3}&
\textbf{2,5,10}&
\textbf{407}\tabularnewline
4&
2,2,5,10&
414\tabularnewline
5&
2,5,5,10,12&
416\tabularnewline
6&
2,5,6,10,10,12&
424\tabularnewline
\hline
\end{tabular}\end{center}

\caption{\label{cap:Table-BIC-1}Selection of the number of groups using BIC
with the model $[a_{i}b_{i}Q_{i}d_{i}]$ of HDDC: data are made of
3 groups with intrinsic dimensions $d_{i}=\{2,5,10\}$.}
\end{table}

Here, we are interested in the selection of the number of groups and
of the intrinsic dimension of the clusters. In this experiment, 
3 Gaussian densities are simulated in $\mathbb{R}^{100}$ according to
the model $[a_{i}b_{i}Q_{i}d_{i}]$ with the following parameters:
$\{ d_{1},d_{2},d_{3}\}=\{2,5,10\}$, $\{\pi_{1},\pi_{2},\pi_{3}\}=\{0.4,0.3,0.3\}$,
$\{ a_{1},a_{2},a_{3}\}=\{150,100,75\}$, $\{ b_{1},b_{2},b_{3}\}=\{15,15,15\}$,
close means and random matrices $Q_{i}$. The dataset was made of
1000 points. Table~\ref{cap:Table-BIC-1} presents the choices of
group intrinsic dimensions for the different values of $k$ and the
corresponding BIC values. First of all, it appears that the criterion
BIC can be successfully used for choosing the number of clusters as
in standard Gaussian mixture models. Indeed, the BIC value
associated to the model $[a_{i}b_{i}Q_{i}d_{i}]$ are computed for different values
of $k$, the number of groups, and BIC indicates that the most likely
value is $k=3$ which is correct. In addition, the intrinsic dimensions
$d_{i}$, estimated by HDDC for $k=3$, are indeed the ones of the
simulated data. It is also interesting to observe the evolution of
the estimation of dimensions $d_{i}$ according to the number of clusters.
For instance, if we consider the case of a mixture of 2 Gaussian densities,
HDDC seems to correctly fit the first 2-dimensional cluster and create
a second cluster made of the two other real groups. In addition, the
estimated dimension of this second cluster is approximately the sum
of the intrinsic dimensions of the two real groups. Similarly, for
$k=4$, HDDC divides the first real group into two new clusters with
intrinsic dimensions equal to $2$. As a conclusion, our approach
for dimension estimation allows to correctly identify
the cluster subspaces.

\subsection{Simulation study: influence of the dimensionality}

\begin{figure}[p]
\begin{center}\includegraphics[%
  width=0.8\columnwidth]{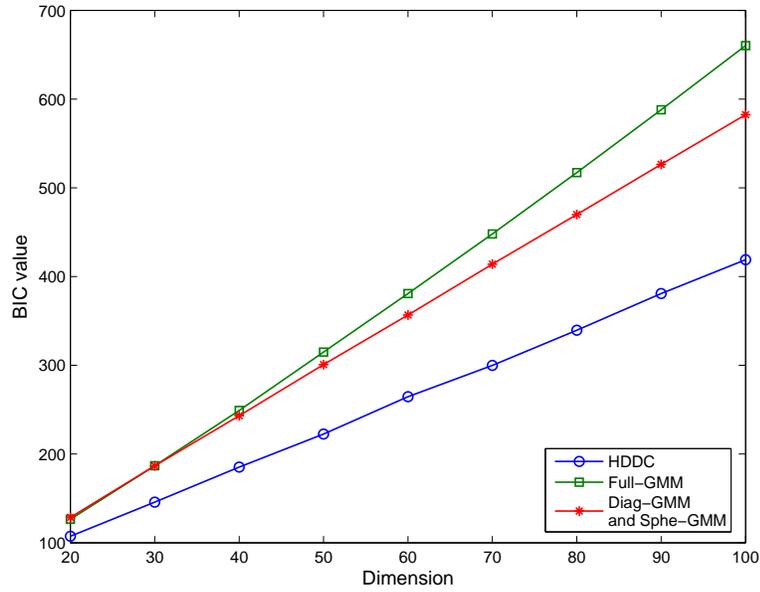}\vspace{-4ex}\end{center}

\caption{\label{cap:Fig-dim-1}Influence of the dimensionality on the BIC
value for different Gaussian mixture models.}
\end{figure}
\begin{figure}[p]
\begin{center}\includegraphics[%
  width=0.8\columnwidth]{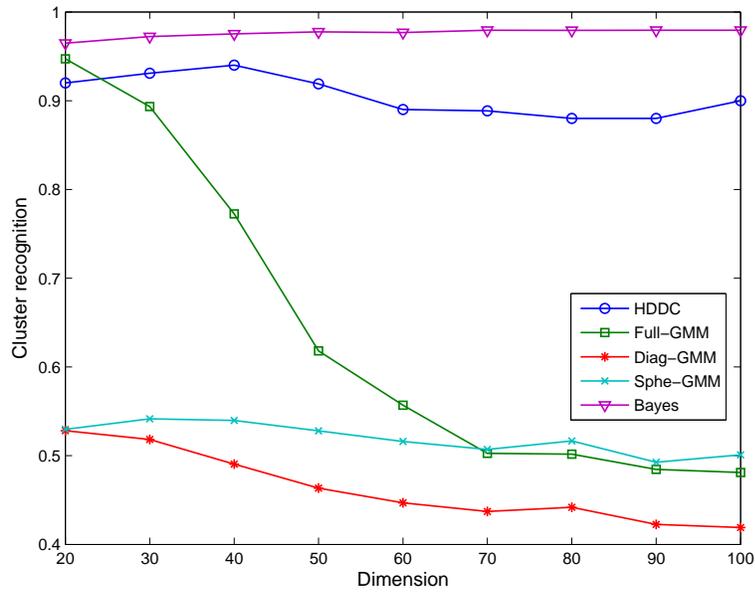}\vspace{-4ex}\end{center}

\caption{\label{cap:Fig-dim-2}Influence of the dimensionality on cluster
recognition rate for different Gaussian mixture models.}
\end{figure}

In this paragraph, we highlight the dimensionality effect on the different
clustering methods. Three
Gaussian densities are simulated in $\mathbb{R}^{p}$, $p=20,...,100$, according
to the model $[a_{i}b_{i}Q_{i}d_{i}]$ with the same parameters as
in the previous experiment. The performance of methods is measured
by the average cluster recognition rate computed on 50 replications.
The studied clustering methods were initialized using the same random
partition. Figures~\ref{cap:Fig-dim-1} and~\ref{cap:Fig-dim-2}
respectively show the influence of the dimensionality on the BIC value
and the cluster recognition rate for different Gaussian mixture models:
model $[a_{i}b_{i}Q_{i}d_{i}]$ of HDDC, Full-GMM, Diag-GMM and Sphe-GMM.
It is not surprising to observe on Figure~\ref{cap:Fig-dim-1} that
BIC selects the model $[a_{i}b_{i}Q_{i}d_{i}]$ as the best model
since the data are simulated according to this model. However, it
interesting to remark that, the more the dimension increases, the
larger the difference between the BIC values of the different models
is, and that in favor of the model $[a_{i}b_{i}Q_{i}d_{i}]$. Figure~\ref{cap:Fig-dim-2}
shows that data dimension does not influence the performance
of HDDC which is very close to the performance of the Bayes decision
rule (computed with the true densities). In addition, HDDC provides
a cluster recognition rate similar to Full-GMM in low dimensions.
Full-GMM is known to be very sensitive to the data dimension
and, indeed, gives bad results as soon as the dimension increases.
The models Diag-GMM and Sphe-GMM cannot correctly fit the data since
they are too parsimonious for this complex dataset. However, one can
observe that Sphe-GMM is not sensitive to the data dimension
whereas Diag-GMM is. To summarize, HDDC is not sensitive to the dimension
and works very well both in low and in high-dimensional spaces. In
addition, the model $[a_{i}b_{i}Q_{i}d_{i}]$ outperforms models requiring
a higher number of parameters (Full-GMM) and models requiring a smaller
number of parameters (Diag-GMM and Sphe-GMM).

\subsection{Simulation study: full rank Gaussian model}

% %
% \begin{figure}
% \begin{center}\includegraphics[%
%   width=1\columnwidth]{images/CSDA_eigenscree.eps}\end{center}
% \caption{\label{cap:FullRank-1}Eigenvalue scree of the simulated covariance
% matrices.}
% \end{figure}

%
\begin{figure}[p]
\begin{center}\includegraphics[%
  width=0.7\columnwidth]{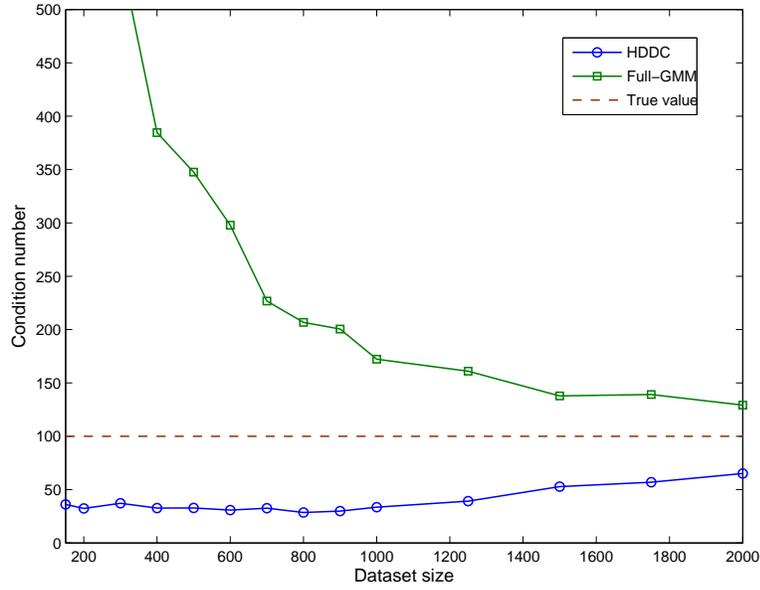}\end{center}
\caption{\label{cap:FullRank-3}Influence of the dataset
size on the condition number for the full rank data.}
\vspace{-4ex}
\end{figure}
\begin{figure}[p]
\begin{center}\includegraphics[%
  width=0.7\columnwidth]{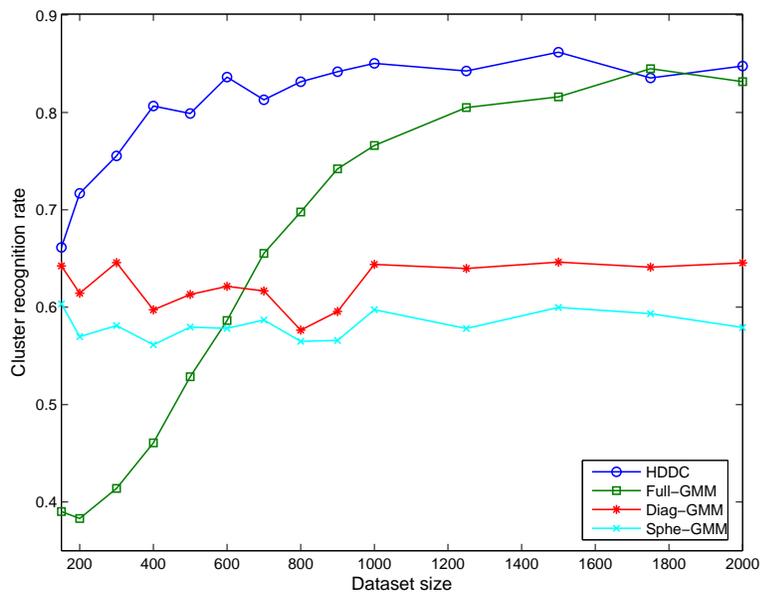}\end{center}
\caption{\label{cap:FullRank-2}Influence of the dataset
size on the cluster recognition rate for the full rank data.}
\vspace{-4ex}
\end{figure}

In this last simulation study, the capacity of HDDC
models to deal with full rank Gaussian data is investigated.
Three Gaussian densities in $\mathbb{R}^{p}$, $p=50$, are simulated with
full rank covariance matrices, \emph{i.e.} according to the model
Full-GMM. The covariance matrices of the groups were different (different
orientations and eigenvalues) but with the same condition number 
fixed to $100$. Recall that the condition number of a matrix is the ratio of its
largest and smallest eigenvalues. 
For this experiment, we used HDDC with the model
$[a_{ij}b_iQ_id_i]$ and the clustering methods associated to the classical
Gaussian models Full-GMM, Diag-GMM and Sphe-GMM. 
In order to observe the
behavior of the studied clustering methods with respect to
the \emph{curse of the dimensionality}, 
the cluster recognition rate is computed for different dataset sizes $n$
since this phenomenon occurs
when the size of the dataset becomes too small compared to the dimension.
As an illustration, Figure~\ref{cap:FullRank-3} presents a comparison 
between the condition number of the estimated covariance matrix associated to the
first group used by the Full-GMM method and the ratio $\hat a_{11}/\hat b_1$, which is the corresponding condition number of the covariance matrix estimated by HDDC, 
for different sizes of the full rank dataset $n=150,...,2000$. 
It appears that, for small datasets (\emph{i.e.} $n$ smaller
than 1000), the condition number of the empirical covariance matrix
(associated to the model Full-GMM) explodes, whereas the condition number
associated to the estimated covariance matrix in the model $[a_{ij}b_iQ_id_i]$
remains stable. 
Figure~\ref{cap:FullRank-2} shows the consequence on the behavior of the studied clustering methods.
First, observe that both Diag-GMM and Sphe-GMM models do not obtain
satisfying results for any dataset size. 
This is due to the fact that the assumptions made
by those models are wrong for the simulated data and they are thus not
able to correctly fit these data. Second, HDDC obtains a similar
cluster recognition rate to the model Full-GMM, which is the model used to
simulate the data, when the dataset size is large (\emph{i.e.} $n$ larger than
1500). Furthermore, HDDC appears to be more efficient to cluster these data than
the model Full-GMM when the dataset size becomes small. Indeed, the cluster
recognition rate of HDDC is almost constant for dataset sizes between 1500 and
500. However, when the dataset size is smaller than 500, the HDDC performance
decreases to the results obtained by the parsimonious
models Diag-GMM and Sphe-GMM. 
These experiments demonstrate that, even with data which are not
favorable to our
model, HDDC is more efficient than both complex and parsimonious models on 
small datasets.

\subsection{Real data study: comparison with variable selection}

\begin{table}[p]
\begin{center}\begin{tabular}{c|c|c}
Model&
Variables&
Cluster recognition rate\tabularnewline
\hline 
Sphe-GMM&
Original&
0.605\tabularnewline
VS-GMM&
Original&
0.925\tabularnewline
Sphe-GMM&
Princ. comp.&
0.605\tabularnewline
VS-GMM&
Princ. comp.&
0.935\tabularnewline
HDDC $[a_{i}b_{i}Q_{i}d_{i}]$&
Original&
\textbf{0.950}\tabularnewline
\end{tabular}\end{center}

\caption{\label{cap:Crabs-data-1}Classification results for the Crabs data:
comparison of different model-based clustering methods.}
\end{table}

\begin{figure}[p]
\begin{center}\includegraphics[%
  width=1\columnwidth]{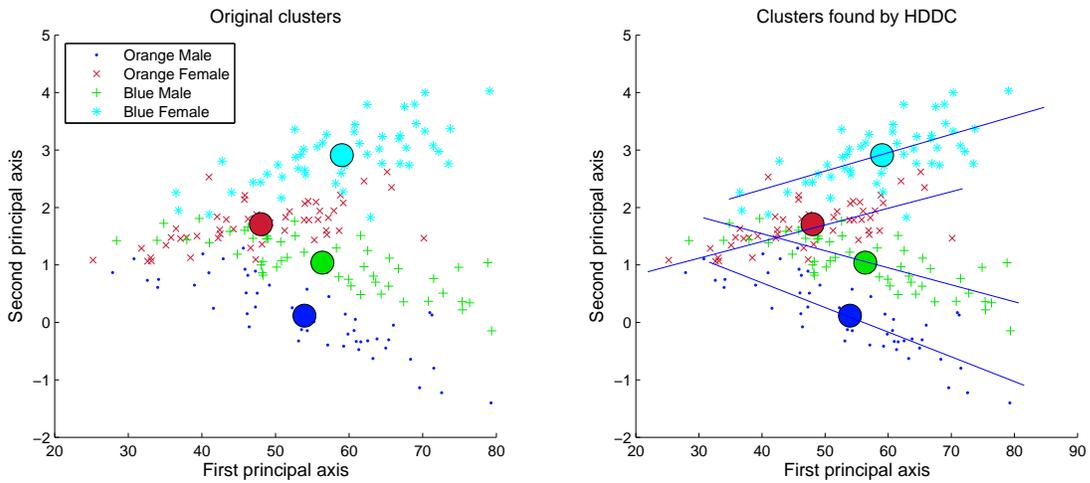}\end{center}
\caption{\label{cap:Crabs-data-2}Clustering results using HDDC: on the left
panel, crabs data projected on the two first principal axes and, on
the right panel, clustering result obtained with the model~$[a_{i}b_{i}Q_{i}d_{i}]$
of HDDC and the estimated specific subspaces of the mixture components
(blue lines).}
\end{figure}

In this experiment, HDDC is compared with the variable selection method
for model-based clustering introduced in~\cite{Raftery05},
and denoted by VS-GMM in the following. The authors considered the variable
selection problem as a model selection problem. Selection is made
using approximate Bayes factors and combined with a greedy search
algorithm. In addition, it is possible to perform this variable selection
on the original variables, but also on the principal components using
PCA as a pre-processing step. In order to compare HDDC to this variable
selection technique, we used the same dataset as in~\cite{Raftery05}.
The Leptograpsus crabs dataset consists of 200 subjects equally distributed
into 4 classes: Orange Male, Orange Female, Blue Male and Blue Female.
There are 5 variables for each subject: width of frontal lip (FL),
rear width (RW), length along the mid-line of the carapace (CL), maximum
of the width of the carapace (CW) and body depth (BD) in mm. The left
panel of Figure~\ref{cap:Crabs-data-2} shows the Crabs data projected
on the two first principal axes and the big circles represent the
cluster means.

%
% \begin{figure}[p]
% \begin{center}\includegraphics[%
 %  bb=8bp 95bp 602bp 757bp,
%   clip,
 %  width=1\columnwidth]{images/fig_crabs_hddc_3.eps}\end{center}

% \caption{\label{cap:Crabs-data-3}The steps of the EM-based algorithm HDDC
% on the Crabs dataset and the estimated specific subspaces of the mixture
% components (blue lines).}
% \end{figure}

Table~\ref{cap:Crabs-data-1} gives the classification error rate
for the classical model Sphe-GMM, the VS-GMM method and HDDC. The
second column of this table indicates on which variables is performed
the clustering. HDDC obtains a cluster recognition rate equal to 95\%
and the variable selection method of Raftery \emph{et al.} obtains
93.5\% whereas the classical model Sphe-GMM obtains a cluster recognition
rate equal to 60.5\%. HDDC found that each cluster lives in a 1-dimensional
subspace embedded into the original 5-dimensional space. The right
panel of Figure~\ref{cap:Crabs-data-2} shows the specific subspaces
(blue lines) of the 4 mixture components obtained with the model~$[a_{i}b_{i}Q_{i}d_{i}]$
of HDDC. For this illustration, the data is projected on the
two first principal components since results obtained with VS-GMM
on these variables are better than on the original ones. It can be observed
that the specific axes of the different clusters are very correlated
and this explains that HDDC provides a better clustering result than
the variable selection method VS-GMM. 
%Figure~\ref{cap:Crabs-data-3}
%presents the 12 steps of the EM algorithm in order to show the evolution
%of the group-specific subspaces into HDDC. 
%We can conclude that HDDC
%is able to fit the clusters in their specific subspaces whereas 
%other methods only select the best dimensions for clustering.
% SANS DOUTE UN PEU TROP GENERAL

\subsection{Real data study: Martian surface characterization }

\begin{figure}[p]
\begin{center}\begin{tabular}{cc}
\includegraphics[%
  width=0.45\columnwidth]{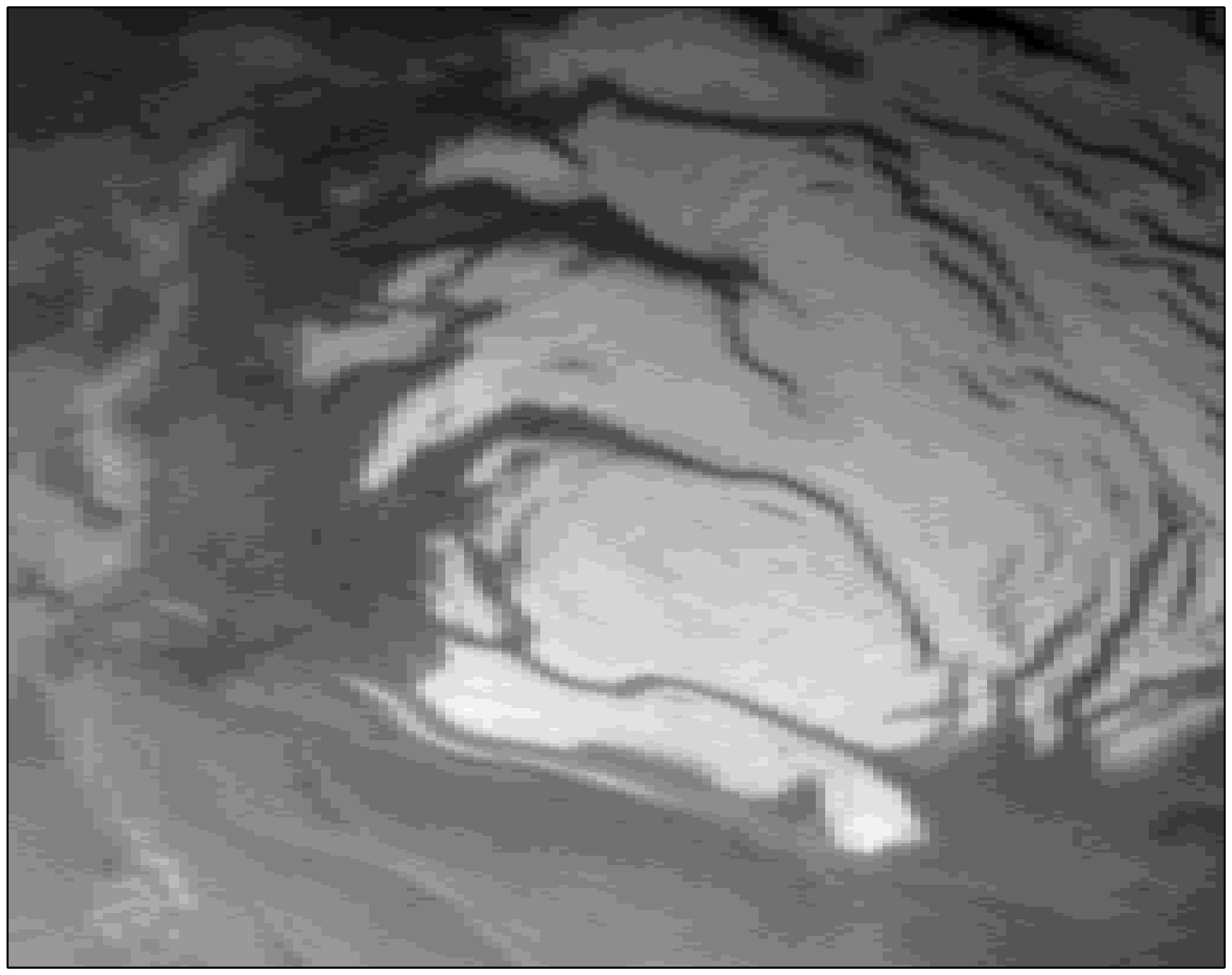}&
\includegraphics[%
  width=0.45\columnwidth]{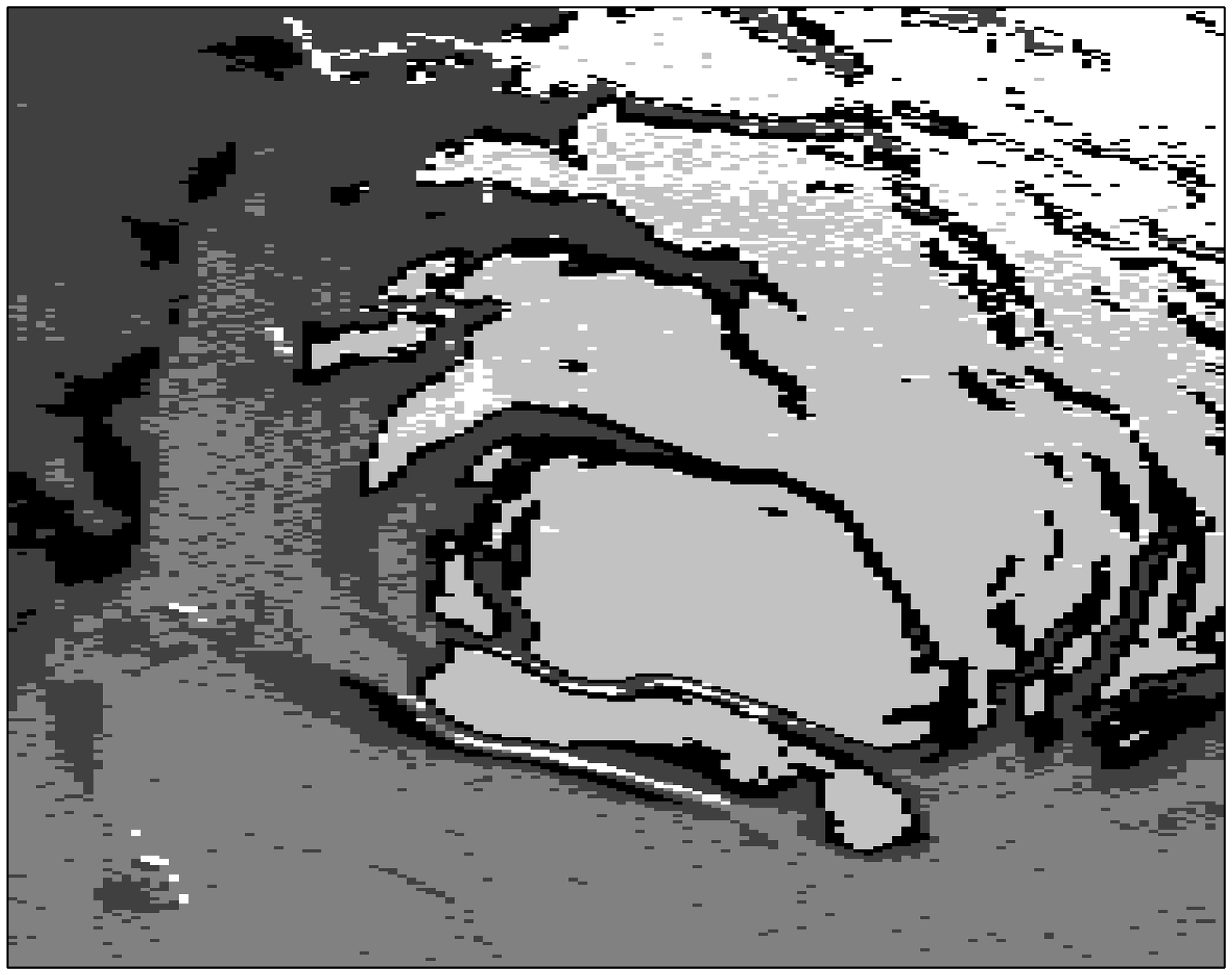}\tabularnewline
\end{tabular}\vspace{-2ex}\end{center}

\caption{\label{cap:Fig-real-data-1}Characterization of the Martian surface
composition using HDDC: on the left, image of the studied zone and,
on the right, segmentation using HDDC on the 256-dimensional spectral
data associated to the image.}
\end{figure}

\begin{figure}[p]
\begin{center}\includegraphics[%
  bb=-62bp 215bp 684bp 564bp,
  clip,
  width=0.9\columnwidth]{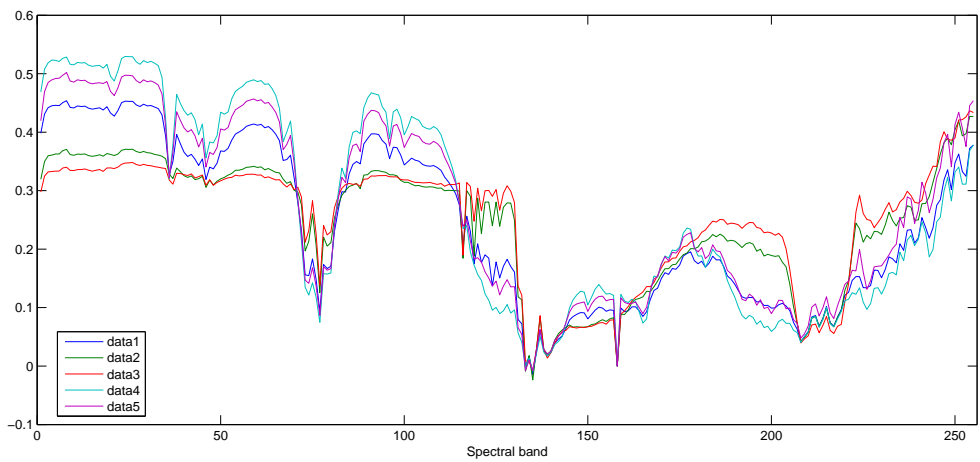}\vspace{-4ex}\end{center}

\caption{\label{cap:Fig-real-data-2}Spectral means of the 5 mineralogical
classes found using HDDC.}
\end{figure}

Here, we propose to use HDDC to analyze and segment images of the
Martian surface. Visible and near infrared imaging spectroscopy is
a key remote sensing technique to study and monitor the system of
the planets. Imaging spectrometers, which are inboard of an increasing
number of satellites, provide high-dimensional hyper-spectral images.
%Constant technological improvements promote the acquisition of dramatically
%expanding image collections.
In March 2004, the OMEGA instrument (Mars
Express, ESA)~\cite{Bibring2004} has collected 310 Gbytes of raw
images. The OMEGA imaging spectrometer has mapped the Martian surface
with a spatial resolution varying between 300 to 3000 meters depending
on spacecraft altitude. It acquires for each resolved pixel the spectrum
from 0.36 to 5.2 \ensuremath{µ}m in 256 contiguous spectral channels.
OMEGA is designed to characterize the composition of surface materials,
discriminating between various classes of silicates, hydrated minerals,
oxides and carbonates, organic frosts and ices. For this experiment,
a $300\times128$ image of the Martian surface is considered and a
256-dimensional spectral observation is associated to each of the
38~400 pixels. The image of the studied zone is presented on the
left panel of Figure~\ref{cap:Fig-real-data-1}. According to the
experts, there are $k=5$ mineralogical classes to identify. 

The right image of Figure~\ref{cap:Fig-real-data-1} shows the segmentation
obtained with the model $[a_{i}b_{i}Q_{i}d_{i}]$ of HDDC. First of
all, observe that the segmentation of HDDC is very precise
on the main part of the image. The poor results of the top right
part of the image are due to the planet curvature and could be corrected.
In particular, the experts of the domain appreciated that our method
is able to detect a mixture of ice and carbonate around
the ice zones (clear zones of the image). Figure~\ref{cap:Fig-real-data-2}
shows the spectral means of the 5 classes and this allows the experts
to determine the mineralogical and molecular composition of each class.
Remind that this study is done without taking into account the
spatial relations between the pixels of a image. A natural extension
of this work is therefore to combine HDDC with the modeling of the
spatial relations using, for example, hidden Markov random fields.
This experiment demonstrates that HDDC can be efficiently used on
real high-dimensional data and with large datasets. In addition, a
main interest of HDDC for this application is 
to provide posterior probabilities that each pixel belongs to the
classes.

\section{Conclusion}

In this paper, new Gaussian mixture models designed
for high-dimensional data are introduced. It is assumed
that the intrinsic dimension
of each mixture component is much smaller than the one
of the original space. In addition, 
outside the specific subspace of each group, the noise variance
is modeled by a single parameter. Additional constraints
can be imposed on the parameters within or between the groups
in order to obtain further regularized models. This parameterization
in the eigenspaces of the mixture components gives rise to an EM-based
clustering method, called High-Dimensional Data Clustering (HDDC).
Experiments on artificial and real datasets demonstrated the effectiveness
of the different models of HDDC compared to classical Gaussian mixture
models. In particular, the model $[a_{i}b_{i}Q_{i}d_{i}]$ provides
very satisfying results for many types of data.

\section*{Acknowledgments}
%The authors are indebted to the referees and the editors for
%their helpful comments and suggestions. They have contributed to
%a greatly improved presentation of the results of this paper.

This work was partially supported by the French department of Research through the project MoViStaR of the \emph{ACI Masse de donn\'ees}.

\appendix

\section{Appendix: parameters estimation}

First of all, we introduce the following useful formulation of the
log-likelihood:\begin{eqnarray}
-2\log(L) & = & \sum_{i=1}^{k}n_{i}\sum_{j=1}^{p}\left(\log(\delta_{ij})+\frac{1}{\delta_{ij}}q_{ij}^{t}W_{i}q_{ij}\right)+c^{st},\label{eq:Likelihood-Flury}\end{eqnarray}
where $\delta_{ij}$ is the $j$th diagonal coefficient of $\Delta_{i}$
and $q_{ij}$ is the $j$th column of $Q_{i}$. We refer to~\cite{Flury84}
for a demonstration of this result.

\subsection{Models with free orientations}

\paragraph*{\noindent Subspace $\mathbb{E}_{i}$:}

\noindent The log-likelihood is to be maximized under the constraint
$q_{ij}^{t}q_{ij}=1$, which is equivalent to finding a saddle point
of the Lagrange function:\[
\mathcal{L}=-2\log(L)-\sum_{j=1}^{p}\theta_{ij}(q_{ij}^{t}q_{ij}-1),\]
where $\theta_{ij}$ are the Lagrange multipliers. Using the expression
(\ref{eq:Likelihood-Flury}) of the log-likelihood, the gradient of
$\mathcal{L}$ with respect to $q_{ij}$ is:\[
\nabla_{q_{ij}}\mathcal{L}=2\frac{n_{i}}{\delta_{ij}}W_{i}q_{ij}-2\theta_{ij}q_{ij},\]
and by multiplying this quantity on the left by $q_{ij}^{t}$, we
obtain:\[
q_{ij}^{t}\nabla_{q_{ij}}\mathcal{L}=0\Leftrightarrow\theta_{ij}=\frac{n_{i}}{\delta_{ij}}q_{ij}^{t}W_{i}q_{ij}.\]
Consequently, $W_{i}q_{ij}=\frac{\theta_{ij}\delta_{ij}}{n_{i}}q_{ij}$
and thus $q_{ij}$ is the eigenvector of $W_{i}$ associated with
the eigenvalue $\lambda_{ij}=\frac{\theta_{ij}\delta_{ij}}{n_{i}}=q_{ij}^{t}W_{i}q_{ij}$.
As the vectors $q_{ij}$ are eigenvectors of the symmetric matrix
$W_{i}$, this implies that $q_{ij}^{t}q_{i\ell}=0$ if $j\neq\ell$.
The log-likelihood can therefore be re-written as follows:\[
-2\log(L)=\sum_{i=1}^{k}n_{i}\left(\sum_{j=1}^{d_{i}}\left(\log(a_{ij})+\frac{\lambda_{ij}}{a_{ij}}\right)+\sum_{j=d_{i}+1}^{p}\left(\log(b_{i})+\frac{\lambda_{ij}}{b_{i}}\right)\right)+c^{st},\]

\noindent and, using the relation $\sum_{j=d_{i}+1}^{p}\lambda_{ij}=\trace(W_{i})-\sum_{j=1}^{d_{i}}\lambda_{ij}$,
we obtain:\begin{equation}
-2\log(L)=\sum_{i=1}^{k}n_{i}\left(\sum_{j=1}^{d_{i}}\log(a_{ij})+(p-d_{i})\log(b_{i})+\frac{\trace(W_{i})}{b_{i}}+\sum_{j=1}^{d_{i}}\left(\frac{1}{a_{ij}}-\frac{1}{b_{i}}\right)\lambda_{ij}\right)+c^{st}.\label{eq:-2log(L)}\end{equation}
Thus, minimizing $-2\log(L)$ with respect to $\lambda_{ij}$ is equivalent
to minimizing the quantity $\sum_{i=1}^{k}n_{i}\sum_{j=1}^{d_{i}}(\frac{1}{a_{ij}}-\frac{1}{b_{i}})\lambda_{ij}$.
Since $(\frac{1}{a_{ij}}-\frac{1}{b_{i}})<0$, $\forall j=1,...,d_{i}$,
$\lambda_{ij}$ must therefore be as larger as possible. Thus, the
column vector $q_{ij}$, $\forall j=1,...,d_{i}$, is estimated by
the eigenvector associated to the $j$th largest eigenvalue of $W_{i}$.

\paragraph*{\noindent Model $[a_{ij}b_{i}Q_{i}d_{i}]$: }

\noindent starting from equation~(\ref{eq:-2log(L)}), the partial
derivative of $-2\log(L)$ with respect to $a_{ij}$ and $b_{i}$
are:\begin{eqnarray*}
-2\frac{\partial\log(L)}{\partial a_{ij}}=n_{i}\left(\frac{1}{a_{ij}}-\frac{\lambda_{ij}}{a_{ij}^{2}}\right) & \textnormal{and} & -2\frac{\partial\log(L)}{\partial b_{i}}=\frac{n_{i}(p-d_{i})}{b_{i}}-\frac{n_{i}}{b_{i}^{2}}\left(\trace(W_{i})-\sum_{j=1}^{d_{i}}\lambda_{ij}\right).\end{eqnarray*}
The condition $\frac{\partial\log(L)}{\partial a_{ij}}=0$ implies
that $\hat{a}_{ij}=\lambda_{ij}$ and the condition $\frac{\partial\log(L)}{\partial b_{i}}=0$
implies that:\[
\hat{b}_{i}=\frac{1}{(p-d_{i})}\left(\trace(W_{i})-\sum_{j=1}^{d_{i}}\lambda_{ij}\right).\]

\paragraph*{\noindent Model $[a_{ij}bQ_{i}d_{i}]$: }

\noindent the partial derivative of $-2\log(L)$ with respect to $b$
is:\[
-2\frac{\partial\log(L)}{\partial b}=\frac{n(p-\xi)}{b}-\frac{1}{b^{2}}\sum_{i=1}^{k}n_{i}\left(\trace(W_{i})-\sum_{j=1}^{d_{i}}\lambda_{ij}\right),\]
and the condition $\frac{\partial\log(L)}{\partial b}=0$ proves that:\[
\hat{b}=\frac{1}{(p-\xi)}\left(\trace(W)-\sum_{i=1}^{k}\hat{\pi}_{i}\sum_{j=1}^{d_{i}}\lambda_{ij}\right).\]

\paragraph*{\noindent Model $[a_{i}b_{i}Q_{i}d_{i}]$: }

\noindent from (\ref{eq:-2log(L)}), the partial derivative of $-2\log(L)$
with respect to $a_{i}$ is:\[
-2\frac{\partial\log(L)}{\partial a_{i}}=\frac{n_{i}d_{i}}{a_{i}}-\frac{n_{i}}{a_{i}^{2}}\sum_{j=1}^{d_{i}}\lambda_{ij},\]
and the condition $\frac{\partial\log(L)}{\partial a_{i}}=0$ implies
that:\[
\hat{a}_{i}=\frac{1}{d_{i}}\sum_{j=1}^{d_{i}}\lambda_{ij}.\]

\paragraph*{\noindent Model $[ab_{i}Q_{i}d_{i}]$:}

\noindent the partial derivative of $-2\log(L)$ with respect to $a$
is:\[
-2\frac{\partial\log(L)}{\partial a}=\frac{n\xi}{a}-\frac{1}{a^{2}}\sum_{i=1}^{k}n_{i}\sum_{j=1}^{d_{i}}\lambda_{ij},\]
 and the condition $\frac{\partial\log(L)}{\partial a}=0$ gives:\[
\hat{a}=\frac{1}{\xi}\sum_{i=1}^{k}\hat{\pi}_{i}\sum_{j=1}^{d_{i}}\lambda_{ij}.\]

\paragraph*{\noindent Model $[a_{j}b_{i}Q_{i}d]$: }

\noindent the partial derivative of $-2\log(L)$ with respect to $a_{j}$
is:\[
-2\frac{\partial\log(L)}{\partial a_{j}}=\frac{n}{a_{j}}-\frac{1}{a_{j}^{2}}\sum_{i=1}^{k}n_{i}\lambda_{ij}.\]
The condition $\frac{\partial\log(L)}{\partial a_{j}}=0$ and the
relation $\sum_{i=1}^{k}n_{i}\lambda_{ij}=n\lambda_{j}$ imply that
$\hat{a}_{j}=\lambda_{j}$.

\subsection{Models with common orientations}

\paragraph*{Subspace $\mathbb{E}_{i}$:}

\noindent Starting from the likelihood expression~(\ref{eq:Likelihood-Flury}),
we can write:\begin{eqnarray*}
-2\log(L) & = & \sum_{i=1}^{k}n_{i}\sum_{j=1}^{d}\left(\log(a_{i})+\frac{1}{a_{i}}q_{j}^{t}W_{i}q_{j}\right)+\sum_{i=1}^{k}n_{i}\sum_{j=d+1}^{p}\left(\log(b_{i})+\frac{1}{b_{i}}q_{j}^{t}W_{i}q_{j}\right)+c^{st},\\
 & = & \sum_{i=1}^{k}n_{i}\left(d\log(a_{i})+(p-d)\log(b_{i})\right)+\sum_{j=1}^{d}q_{j}^{t}Aq_{j}+\sum_{j=d+1}^{p}q_{j}^{t}Bq_{j}+c^{st},\end{eqnarray*}
where $A=\sum_{i=1}^{k}\frac{n_{i}}{a_{i}}W_{i}$ and $B=\sum_{i=1}^{k}\frac{n_{i}}{b_{i}}W_{i}$.
Note that $\sum_{j=d+1}^{p}q_{j}^{t}Bq_{j}$ can be written using
the trace of $B$: $\sum_{j=d+1}^{p}q_{j}^{t}Bq_{j}=\trace(B)-\sum_{j=1}^{d}q_{j}^{t}Bq_{j}$.
This yields:\begin{eqnarray}
-2\log(L) & = & \sum_{i=1}^{k}n_{i}\left(d\log(a_{i})+(p-d)\log(b_{i})\right)-\sum_{j=1}^{d}q_{j}^{t}(B-A)q_{j}+\trace(B)+c^{st}.\label{eq:-2log(L)_common_Q}\end{eqnarray}
Consequently, the gradient of $\mathcal{L}=-2\log(L)-\sum_{j=1}^{p}\theta_{j}(q_{j}^{t}q_{j}-1)$
with respect to $q_{j}$ is: \[
\nabla_{q_{j}}\mathcal{L}=-2(B-A)q_{j}-2\theta_{j}q_{j},\]
where $\theta_{j}$ is the $j$th Lagrange multiplier. The relation
$\nabla_{q_{j}}\mathcal{L}=0$ is equivalent to $(B-A)q_{j}=-\theta_{j}q_{j}$
which means that $q_{j}$ is eigenvector of the matrix $(B-A)$. In
order to minimize the quantity $-2\log(L)$, the $d$ first columns
of $Q$ must be the eigenvectors associated with the $d$ largest
eigenvalues of $(B-A)$.

\paragraph*{\noindent Model $[a_{i}b_{i}Qd]$: }

\noindent Starting from equation~(\ref{eq:-2log(L)_common_Q}), the
partial derivatives of $-2\log(L)$ with respect to $a_{i}$ and $b_{i}$
are:\[
-2\frac{\partial\log(L)}{\partial a_{i}}=\frac{n_{i}d}{a_{i}}-\frac{n_{i}}{a_{i}^{2}}\sum_{j=1}^{d}q_{j}^{t}W_{i}q_{j}\,\,\textnormal{and}\,-2\frac{\partial\log(L)}{\partial b_{i}}=\frac{n_{i}(p-d)}{b_{i}}-\frac{n_{i}}{b_{i}^{2}}\left(\trace(W_{i})-\sum_{j=1}^{d}q_{j}^{t}W_{i}q_{j}\right).\]
The condition $\frac{\partial\log(L)}{\partial a_{i}}=0$ and $\frac{\partial\log(L)}{\partial b_{i}}=0$
give respectively:\begin{eqnarray*}
\hat{a}_{i}(Q)=\frac{1}{d}\sum_{j=1}^{d}q_{j}^{t}W_{i}q_{j} & \textnormal{and} & \hat{b}_{i}(Q)=\frac{1}{(p-d)}\left(\trace(W_{i})-\sum_{j=1}^{d}q_{j}^{t}W_{i}q_{j}\right).\end{eqnarray*}

\paragraph*{\noindent Model $[a_{i}bQd]$: }

\noindent The partial derivative of $-2\log(L)$ with respect to $b$
is:\[
-2\frac{\partial\log(L)}{\partial b}=\frac{n(p-d)}{b}-\frac{n}{b^{2}}\left(\trace(W)-\sum_{j=1}^{d}q_{j}^{t}Wq_{j}\right),\]
and the condition $\frac{\partial\log(L)}{\partial b}=0$ implies
that:\[
\hat{b}(Q)=\frac{1}{(p-d)}\left(\trace(W)-\sum_{j=1}^{d}q_{j}^{t}Wq_{j}\right).\]

\paragraph*{\noindent Model $[ab_{i}Qd]$: }

\noindent The partial derivative of $-2\log(L)$ with respect to $a$
is:\[
-2\frac{\partial\log(L)}{\partial a}=\frac{nd}{a}-\frac{n}{a^{2}}\sum_{j=1}^{d}q_{j}^{t}Wq_{j},\]
and the condition $\frac{\partial\log(L)}{\partial a}=0$ proves that:\[
\hat{a}(Q)=\frac{1}{d}\sum_{j=1}^{d}q_{j}^{t}Wq_{j}.\]

\subsection{Models with common covariance matrices}

\paragraph*{\noindent Subspace $\mathbb{E}_{i}$:}

\noindent The log-likelihood can be written as follows:\[
-2\log(L)=n\left(\sum_{j=1}^{d}\log(a_{j})+(p-d)\log(b)+\frac{\trace(W)}{b}+\sum_{j=1}^{d}\left(\frac{1}{a_{j}}-\frac{1}{b}\right)q_{j}^{t}Wq_{j}\right)+c^{st}.\]
 The gradient of $\mathcal{L}=-2\log(L)-\sum_{j=1}^{p}\theta_{j}(q_{j}^{t}q_{j}-1)$
with respect to $q_{j}$ is: \[
\nabla_{q_{j}}\mathcal{L}=2n(\frac{1}{a_{j}}-\frac{1}{b})Wq_{j}-2\theta_{j}q_{j},\]
where $\theta_{j}$ is the $j$th Lagrange multiplier. The relation
$\nabla_{q_{j}}\mathcal{L}=0$ implies that $q_{j}$is eigenvector
of $W$. In order to minimize $-2\log(L)$, the first columns of $Q$
must be the eigenvectors associated to the $d$ largest eigenvalues
of $W$.

\paragraph*{\noindent Model $[a_{j}bQd]$: }

\noindent The partial derivatives of $-2\log(L)$ with respect to
$a_{j}$ and $b$ are:\begin{eqnarray*}
-2\frac{\partial\log(L)}{\partial a_{j}}=\frac{n}{a_{j}}-\frac{n}{a_{j}^{2}}q_{j}^{t}Wq_{j} & \textnormal{and} & -2\frac{\partial\log(L)}{\partial b}=\frac{n(p-d)}{b}-\frac{n}{b^{2}}\sum_{j=d+1}^{p}q_{j}^{t}Wq_{j}.\end{eqnarray*}
The condition $\frac{\partial\log(L)}{\partial a_{i}}=0$ implies
that $\hat{a}_{j}=\lambda_{j}$. The combination of the condition
$\frac{\partial\log(L)}{\partial b}=0$ with the relation $\sum_{j=d+1}^{p}\lambda_{j}=\trace(W)-\sum_{j=1}^{d}\lambda_{j}$
gives the estimator of $b$:\[
\hat{b}=\frac{1}{(p-d)}\left(\trace(W)-\sum_{j=1}^{d}\lambda_{j}\right).\]

\paragraph*{\noindent Model $[abQd]$: }

\noindent The partial derivatives of $-2\log(L)$ with respect to
$a$ is:\[
-2\frac{\partial\log(L)}{\partial a}=\frac{nd}{a}-\frac{n}{a^{2}}\sum_{j=1}^{d}q_{j}^{t}Wq_{j},\]
and the condition $\frac{\partial\log(L)}{\partial a}=0$ implies
that: \[
\hat{a}=\frac{1}{d}\sum_{j=1}^{d}\lambda_{j}.\]

\bibliographystyle{plain}
\bibliography{biblio_CSDA}

\end{document}